\documentclass[a4paper]{article}

\usepackage{a4wide}
\usepackage{latexsym}
\usepackage[UKenglish]{babel}
\usepackage{graphicx,subfigure}
	\graphicspath{{./}{figures/}}
\usepackage[colorlinks=true,linkcolor=blue,citecolor=red]{hyperref}
\usepackage{xcolor}
\usepackage{amsfonts,amsmath,amsthm,bbm,cool}

\newtheorem{theorem}{Theorem}[section]
\newtheorem{proposition}[theorem]{Proposition}
\theoremstyle{remark}\newtheorem{remark}[theorem]{Remark}

\newcommand{\abs}[1]{\left\lvert#1\right\rvert}
\newcommand{\ave}[1]{\overline{#1}}
\newcommand{\bracket}[1]{\left\langle#1\right\rangle}
\newcommand{\be}{\begin{equation}}
\newcommand{\ee}{\end{equation}}
\newcommand{\C}{\mathcal{C}}

\newcommand{\F}{\mathcal{F}}

\newcommand{\R}{\mathbb{R}}
\newcommand{\rangeth}{I_\Theta}
\newcommand{\Var}[1]{\operatorname{Var}(#1)}

\allowdisplaybreaks

\begin{document}
\title{Boltzmann-type models with uncertain binary interactions}

\author{	Andrea Tosin \\
		{\small 	Department of Mathematical Sciences ``G. L. Lagrange''} \\
		{\small Politecnico di Torino, Torino, Italy} \\
		{\small\tt andrea.tosin@polito.it} \\[5mm]
		Mattia Zanella \\
		{\small	Department of Mathematical Sciences ``G. L. Lagrange''} \\
		{\small Politecnico di Torino, Torino, Italy} \\
		{\small\tt mattia.zanella@polito.it}}
\date{}

\maketitle

\begin{abstract}
In this paper we study binary interaction schemes with uncertain parameters for a general class of Boltzmann-type equations with applications in classical gas and aggregation dynamics. We consider deterministic (i.e., \emph{a priori} averaged) and stochastic kinetic models, corresponding to different ways of understanding the role of uncertainty in the system dynamics, and compare some thermodynamic quantities of interest, such as the mean and the energy, which characterise the asymptotic trends. Furthermore, via suitable scaling techniques we derive the corresponding deterministic and stochastic Fokker-Planck equations in order to gain more detailed insights into the respective asymptotic distributions. We also provide numerical evidences of the trends estimated theoretically by resorting to recently introduced structure preserving uncertainty quantification methods.

\medskip

\noindent{\bf Keywords:} Uncertainty quantification, deterministic and stochastic kinetic equations, Boltz\-mann and Fokker-Planck equations, structure preserving schemes \\

\noindent{\bf Mathematics Subject Classification:} 35Q20, 35Q70, 35Q84
\end{abstract}


\section{Introduction}
\label{sect:intro}
In recent years kinetic theory has emerged as one of the most natural theoretical frameworks for the description of a wide range of non-classical phenomena spanning from opinion dynamics~\cite{albi2016CHAPTER,albi2017KRM,cristiani2017PREPRINT,during2009PRSA,during2015PRSA,toscani2006CMS} to socio-economic systems~\cite{cordier2005JSP,tosin2014MPE2013}, swarm dynamics~\cite{carrillo2009KRM,carrillo2010SIMA,ha2008KRM}, crowd dynamics~\cite{agnelli2015M3AS,albi2016SIAP,degond2013JSP,degond2013KRM,festa2017KRM}, vehicular traffic~\cite{fermo2014DCDSS,freguglia2017CMS,gunther2002MCM,puppo2016CMS,puppo2017KRM,visconti2017MMS} and biological systems~\cite{chalub2004MM,chauviere2007NHM,perthame2004AM}. The underlying idea is that systems composed by a sufficiently large number of particles/agents can be described through the laws of statistical physics as it happens in the classical theory of rarefied gases. Such an analogy allows one to exploit the sound theoretical background of the Boltzmann equation both at the analytical and at the numerical level~\cite{cercignani1988BOOK,during2009RMUP,pareschi2013BOOK}.

One of the key aspects of the kinetic framework is its flexibility in the description of complex dynamics at different scales. Vlasov-Fokker-Planck and hydrodynamic equations can be derived, which may allow for the explicit computation of the stationary states and of minimal energy states of the system~\cite{furioli2017M3AS}. The rigorous derivation of mean-field and fluid dynamic models from the Boltzmann equation is a deeply fascinating issue. Without intending to review all the literature on the topic, we indicate~\cite{carrillo2010MSSET,degond2008KRM,degond2013MAA,duan2010CMP,ha2008KRM} as representative works of the hierarchy of scales described by kinetic modelling.

It is well known that models of real world phenomena are usually affected by ineradicable uncertainties in some of their relevant parameters, which often are known only statistically. In collective phenomena, such as those encountered in the aforementioned applications, this issue is particularly evident since the dynamics are inferred mostly from empirical observations and models are often calibrated by means of experimental data~\cite{ballerini2008AB,bongini2017M3AS,katz2011PNAS}. In the context of kinetic modelling, this issue may be translated in a general uncertainty in the binary interaction rules, which may include a random variable $\theta$, taking values in a set $\rangeth\subseteq\R$, with known probability distribution $\rho=\rho(\theta):\rangeth\to\R_+$. For the sake of simplicity, we will assume $\theta$ scalar but the same arguments may be repeated for a vector-valued $\theta$ modelling several sources of uncertainty in the interaction rules.

The uncertainty in the binary interactions affects the evolution of the kinetic distribution function, thus the observable trends of the system. In this paper we examine and compare two possible manners of quantifying the impact of the uncertainty at the collective level, corresponding to two different ways of understanding its role in the system dynamics. On one hand, we consider \emph{$\theta$-averaged} binary interactions and study the evolution of the distribution function $g=g(t,\,v)$, $v$ being the generic microscopic state of the particles, ruled by the Boltzmann-type equation
\begin{subequations}
\begin{equation}
	\partial_tg(t,\,v)=\int_{\rangeth}Q(g,\,g)(t,\,v;\,\theta)\rho(\theta)\,d\theta,
	\label{eq:g.intro}
\end{equation}
where $Q$ is the collisional operator, which depends on $\theta$ for so do the interaction rules, and the integral gives its expected value with respect to the law of $\theta$. Notice that the distribution function $g$ is \emph{deterministic}, because it does not depend on $\theta$. From the modelling point of view this corresponds to thinking of $\theta$ as a quantity which varies in each binary interaction, so that we study the statistical evolution of the system by referring to the representative $\theta$-averaged binary interaction. On the other hand, we consider binary interactions \emph{parametrised} by $\theta$ and study the evolution of the distribution function $f=f(t,\,v;\,\theta)$ obeying the Boltzmann-type equation
\begin{equation}
	\partial_tf(t,\,v;\,\theta)=Q(f,\,f)(t,\,v;\,\theta).
	\label{eq:f.intro}
\end{equation}
\end{subequations}
In this case $f$ is a \emph{stochastic} distribution, because it still depends on the uncertainty intrinsic in $\theta$. From the modelling point of view this amounts to thinking of $\theta$ as a quantity constant in each binary interaction, whose precise value is however unknown. Average collective dynamics can then be obtained by averaging the family of distribution functions $\{f(t,\,v;\,\cdot)\}$ with respect to the law of $\theta$, each $f(t,\,v;\,\theta)$ corresponding to a more or less probable evolution of the system depending on the distribution of $\theta$.

In this paper we will observe how termodinamic moments and large time trends of the solutions to~\eqref{eq:g.intro}-\eqref{eq:f.intro} ($\theta$-averaged in the latter case) may be either indistinguishable or considerably different for the same binary interaction model. In particular, we will exemplify this fact by means of quite different classical models, such as the Kac model for a caricature of gas molecules with uncertain angle of collision~\cite{kac1959BOOK} and a general aggregation/consensus model with uncertain interactions, see~\cite{toscani2006CMS} for its deterministic version. Since the technical complexity of the kinetic equations allows for the explicit computation of analytical solutions only in particular cases, we will also make use of suitable numerical techniques coming from the uncertainty quantification (UQ)~\cite{dimarco2017CHAPTER,hu2017CHAPTER,mishra2013CHAPTER} for quantitatively describing more complex uncertain dynamics, such as opinion dynamics with an uncertain bounded confidence threshold.

In more detail, the paper is organised as follows: in Section~\ref{sect:Boltzmann} we introduce the general uncertain binary dynamics and the corresponding deterministic and stochastic Boltzmann-type equations. Then, in such a framework, we study the Kac model and general aggregation/consensus models under the two different perspectives discussed above. In Section~\ref{sect:Fokker-Planck} we apply the quasi-invariant interaction limit~\cite{toscani2006CMS} to obtain Fokker-Planck-type asymptotic models approximating the large time solutions of the corresponding Boltzmann models in a balanced transport-diffusion regime. Taking advantage of such a tool we then discuss the deterministic and the stochastic stationary distributions produced by the uncertain microscopic dynamics for the inelastic Kac model~\cite{pulvirenti2004JSP} and for general diffusive aggregation/consensus models. In Section~\ref{sect:numerics} we present some numerical examples obtained by means of recently developed structure preserving stochastic collocation schemes. Finally, in Section~\ref{sect:conclusions} we collect some concluding remarks about the contents of the paper.

\section{Boltzmann-type equations}
\label{sect:Boltzmann}
\subsection{General uncertain binary interactions}
Following~\cite{pareschi2013BOOK}, we consider linear binary interaction models on the real line leading from pre-interaction states $v,\,w\in V\subseteq\R$ to post-interaction states $v^\ast,\,w^\ast\in V$ of the form
\begin{equation}
	v^\ast=p_1(\theta)v+q_1(\theta)w, \qquad w^\ast=p_2(\theta)v+q_2(\theta)w,
	\label{eq:binary_general}
\end{equation}
where $p_i,\,q_i$, $i=1,\,2$, are random model parameters depending on a real-valued random variable $\theta\in\rangeth\subseteq\R$. We assume that the latter is distributed according to a known probability density $\rho:\R\to\R_+$, such that $\rho(\theta)\geq 0$ a.e. in $\R$, $\operatorname{supp}\rho\subseteq\rangeth$ and $\int_{\rangeth}\rho(\theta)\,d\theta=1$.

For the moment, in~\eqref{eq:binary_general} we do not include additional stochastic effects leading to diffusive-like contributions at a collective level. We postpone them to the Fokker-Planck asymptotic analysis of Section~\ref{sect:Fokker-Planck}.

As anticipated in the Introduction, in order to deal with the uncertainty contained in~\eqref{eq:binary_general} in a kinetic description of the particle system, essentially two approaches are possible. On one hand, one can average the effect of each binary interaction, which corresponds to considering the evolution of the kinetic distribution function $g=g(t,\,v):\R_+\times V\to\R_+$ ruled by the Boltzmann-type equation (in weak form)
\begin{align}
	\begin{aligned}[b]
		\frac{d}{dt} & \int_V\varphi(v)g(t,\,v)\,dv	\\
		&=\frac{1}{2}\iint_{V^2}\int_{\rangeth}\left(\varphi(v^\ast)+\varphi(w^\ast)
			-\varphi(v)-\varphi(w)\right)\rho(\theta)g(t,\,v)g(t,\,w)\,d\theta\,dv\,dw \\
		&=\frac{1}{2}\iint_{V^2}\left(\ave{\varphi(v^\ast)+\varphi(w^\ast)}
			-\varphi(v)-\varphi(w)\right)g(t,\,v)g(t,\,w)\,dv\,dw,
	\end{aligned}
	\label{eq:Boltzmann_1}
\end{align}
cf.~\cite{pareschi2013BOOK}, where $\varphi:V\to\R$ is a test function and $\ave{(\cdot)}$ denotes, here and henceforth, the average with respect to $\theta$, i.e., $\ave{(\cdot)}:=\int_{\rangeth}(\cdot)\rho(\theta)\,d\theta$. Since the distribution function $g$ does not depend on the random variable $\theta$, we call~\eqref{eq:Boltzmann_1} the \emph{deterministic model} associated to the interaction rules~\eqref{eq:binary_general}.

On the other hand, one can consider a Boltzmann-type equation for the kinetic distribution function $f=f(t,\,v;\,\theta):\R_+\times V\times\rangeth\to\R_+$ \emph{parametrised} by (viz. conditioned to) $\theta$:
\begin{align}
	\begin{aligned}[b]
		\frac{d}{dt} & \int_V\varphi(v)f(t,\,v;\,\theta)\,dv  \\
		&= \frac{1}{2}\iint_{V^2}\left(\varphi(v^\ast)+\varphi(w^\ast)-\varphi(v)-\varphi(w)\right)
			f(t,\,v;\,\theta)f(t,\,w;\,\theta)\,dv\,dw.
	\end{aligned}
	\label{eq:Boltzmann_2}
\end{align}
In this case the uncertainty caused by $\theta$ affects the global dynamics and has to be averaged \emph{a posteriori} directly at the collective level. This is the typical approach followed for the uncertainty quantification. In contrast to the previous case, we call~\eqref{eq:Boltzmann_2} the \emph{stochastic model} associated to the interaction rules~\eqref{eq:binary_general}.

Notice that choosing $\varphi(v)=1$ in~\eqref{eq:Boltzmann_1},~\eqref{eq:Boltzmann_2} yields
\be \label{eq:mass}
\frac{d}{dt}\int_Vg(t,\,v)\,dv=0, \qquad \frac{d}{dt}\int_Vf(t,\,v;\,\theta)\,dv=0, 
\ee
hence both $g(t,\,\cdot)$ and $f(t,\,\cdot;\,\theta)$ can be regarded as probability densities over the microscopic state $v$ for all $t>0$ and, in the latter case, also all $\theta\in\rangeth$ if they are so at $t=0$. The evolution of higher-order thermodynamic-like moments of the distribution functions is however not expected to be the same in general, considering that the moments of $f$ depend explicitly on $\theta$. For instance, for the mean we have
\be \label{eq:momentum}
m_g(t):=\int_Vvg(t,\,v)\,dv, \qquad m_f(t;\,\theta):=\int_Vvf(t,\,v;\,\theta)\,dv 
\ee
and for the second order moment (related to the $v$-variance)
\be \label{eq:energy}
E_g(t):=\int_Vv^2g(t,\,v)\,dv, \qquad E_f(t;\,\theta):=\int_Vv^2f(t,\,v;\,\theta)\,dv. 
\ee
In the case of $f$, we can further consider mean quantities with respect to $\theta$, such as
\begin{align*}
	\bar{f}(t,\,v) &:= \ave{f(t,\,v;\,\theta)} \\
	m_{\bar{f}}(t) &:= \int_Vv\bar{f}(t,\,v)\,dv=\ave{m_f(t;\,\theta)} \\
	E_{\bar{f}}(t) &:= \int_Vv^2\bar{f}(t,\,v)\,dv=\ave{E_f(t;\,\theta)}.
\end{align*}

In the forthcoming sections we will investigate in some detail analogies and differences between the dynamics predicted by~\eqref{eq:Boltzmann_1} and~\eqref{eq:Boltzmann_2}.

\subsection{The Kac model}
\label{sect:Kac}
An example falling in the general framework~\eqref{eq:binary_general} is provided by the celebrated Kac model~\cite{kac1959BOOK}, which describes binary collisions among indistinguishable molecules of a spatially homogeneous gas. In this case, $v,\,w$ are the velocities of the colliding molecules, which for the sake of simplicity are assumed one-dimensional in the state space $V=\R$.

The Kac model is obtained from~\eqref{eq:binary_general} with the choices
\begin{equation}
	p_1(\theta)=q_2(\theta)=\cos{\theta}, \qquad p_2(\theta)=-q_1(\theta)=\sin{\theta},
	\label{eq:Kac}
\end{equation}
$\theta$ being a random parameter uniformly distributed in $\rangeth=[0,\,2\pi]$, i.e. $\theta\sim\mathcal{U}(0,\,2\pi)$ with consequently $\rho(\theta)=\frac{1}{2\pi}\mathbbm{1}_{[0,\,2\pi]}(\theta)$. The interaction rules~\eqref{eq:binary_general} together with~\eqref{eq:Kac} imply $(v^\ast)^2+(w^\ast)^2=v^2+w^2$, i.e. the energy is conserved in each binary collision, hence also globally. The same is instead not true for the mean velocity, therefore it is interesting to study the evolution of $m_g$ and $m_{\bar{f}}$ according to~\eqref{eq:Boltzmann_1},~\eqref{eq:Boltzmann_2}, respectively.

Setting, to fix the ideas, $m_g(0)=m_f(0;\,\theta)=1$ for all $\theta\in [0,\,2\pi]$ and taking $\varphi(v)=v$ in~\eqref{eq:Boltzmann_1},~\eqref{eq:Boltzmann_2} we obtain
$$ m_g(t)=e^{-t}, \qquad m_f(t;\,\theta)=e^{(\cos{\theta}-1)t}, $$
whence in particular
$$ m_{\bar{f}}(t)=\int_{\rangeth}m_f(t;\,\theta)\rho(\theta)\,d\theta=\frac{1}{2\pi}\int_{0}^{2\pi}e^{(\cos{\theta}-1)t}\,d\theta. $$
Since $e^{(\cos{\theta}-1)t}\leq 1$ for all $t\geq 0$ and $e^{(\cos{\theta}-1)t}$ converges pointwise to zero for all $\theta\in (0,\,2\pi)$ when $t\to+\infty$, by dominated convergence we get
$$ \lim_{t\to+\infty}m_{\bar{f}}(t)=\frac{1}{2\pi}\int_{0}^{2\pi}\lim_{t\to+\infty}e^{(\cos{\theta}-1)t}\,d\theta=0, $$
therefore $m_g$ and $m_{\bar{f}}$ have the same asymptotic trend. However, considering that $\cos{\theta}-1\geq-\frac{1}{2}\theta^2$, we further discover
\begin{align*}
	\int_{0}^{2\pi}e^{(\cos{\theta}-1)t}\,d\theta &= 2\int_{0}^{\pi}e^{(\cos{\theta}-1)t}\,d\theta \\
	&\geq 2\int_{0}^{\pi}e^{-\frac{t}{2}\theta^2}\,d\theta=\sqrt{\frac{2\pi}{t}}\Erf{\pi\sqrt{\frac{t}{2}}},
\end{align*}
where $\Erf{x}:=\frac{2}{\sqrt{\pi}}\int_{0}^{x}e^{-y^2}\,dy$ is the error function. Consequently
$$ m_{\bar{f}}(t)\geq\frac{1}{\sqrt{2\pi t}}\Erf{\pi\sqrt{\frac{t}{2}}}=O\left(\frac{1}{\sqrt{t}}\right) \qquad (t\to+\infty), $$
indicating that the convergence to zero of $m_{\bar{f}}$ is by far much slower than that of $m_g$, which is instead exponentially fast.

\subsection{Aggregation/consensus models}
\label{sect:consensus}
Let us now consider in~\eqref{eq:binary_general} the \emph{symmetric} case, i.e. $p_1=q_2=:p$, $p_2=q_1=:q$, with additionally $p(\theta)=1-q(\theta)$. Hence the binary interaction scheme simplifies to 
\begin{equation}
	v^\ast=v+q(\theta)(w-v), \qquad w^\ast=w+q(\theta)(v-w),
	\label{eq:binary_symmetric}
\end{equation}
which can be regarded as a prototype of the interactions used in a wide variety of applications such as, e.g., opinion dynamics models~\cite{albi2016CHAPTER,during2015PRSA,toscani2006CMS}, vehicular traffic models~\cite{gunther2002MCM,herty2010KRM}, wealth distribution models~\cite{cordier2005JSP,slanina2004PRE}, spatially homogeneous swarming models~\cite{albi2013AML,carrillo2010MSSET}. In this case, since $v^\ast+w^\ast=v+w$, the average microscopic state is conserved in each interaction, hence also globally. Therefore, if we assume that the initial distribution is the same for both $g$ and $f$, being in particular independent of $\theta$, and has zero mean then we get immediately $m_g(t)=m_f(t;\,\theta)=m_{\bar{f}}(t)=0$ for all $t>0$ and all $\theta\in\rangeth$.

This makes it possible to study the average asymptotic trends induced by~\eqref{eq:binary_symmetric} by looking only at $E_g$ and $E_{\bar{f}}$, which coincide with the variance of $g$ and $\bar{f}$, respectively. In particular, we say that the rules~\eqref{eq:binary_symmetric} lead to an \emph{aggregation}, or \emph{consensus}, for models~\eqref{eq:Boltzmann_1},~\eqref{eq:Boltzmann_2} if
\begin{equation}
	\lim_{t\to+\infty}W_2(g(t),\,\delta_0)=0, \qquad \lim_{t\to+\infty}W_2(\bar{f}(t),\,\delta_0)=0,
	\label{eq:Wasserstein}
\end{equation}
respectively, where $\delta_0$ is the Dirac distribution centred at $v=0$ and $W_2$ denotes the $2$-Wasserstein distance in the space of probability measures, see e.g.~\cite{ambrosio2008BOOK}. Conditions~\eqref{eq:Wasserstein} are met, in particular, if $E_g(t),\,E_{\bar{f}}(t)\to 0$ when $t\to+\infty$, because $W_2(g(t),\,\delta_0)\leq E_g(t)$ and $W_2(\bar{f}(t),\,\delta_0)\leq E_{\bar{f}}(t)$ for all $t\geq 0$.

Setting $E_g(0)=E_f(0;\,\theta)=1$ for all $\theta\in\rangeth$ and taking $\varphi(v)=v^2$ in~\eqref{eq:Boltzmann_1},~\eqref{eq:Boltzmann_2} we obtain, after standard calculations,
$$ E_g(t)=e^{2\left(\ave{q^2}-\ave{q}\right)t}, \qquad E_f(t;\,\theta)=e^{2q(\theta)(q(\theta)-1)t}. $$
If we fix specifically
\begin{equation}
	q(\theta)=q_0+\lambda\theta, \qquad \lambda>0,
	\label{eq:q}
\end{equation}
and we assume further that the random variable $\theta$ is symmetric, hence such that $\ave{\theta}=0$, with $\Var{\theta}<+\infty$ then we get in particular
\begin{equation}
	E_g(t)=e^{2\left(q_0^2-q_0+\lambda^2\Var{\theta}\right)t}, \qquad
		E_f(t;\,\theta)=e^{2\left(q_0^2-q_0+\lambda^2\theta^2+\lambda(2q_0-1)\theta\right)t}.
	\label{eq:Eg_Ef}
\end{equation}

The expression~\eqref{eq:q} corresponds to assuming that the uncertainty in the coefficient $q$ is due to stochastic fluctuations with variance proportional to $\lambda^2$ around a constant value $q_0$. Since from~\eqref{eq:binary_symmetric}-\eqref{eq:q} without uncertainty ($\lambda=0$) it results $\abs{w^\ast-v^\ast}=\abs{1-2q_0}\cdot\abs{w-v}$, one typically requires $\abs{1-2q_0}<1$ so as to have \emph{contractive} binary interactions fostering aggregation/consensus. Hence in the following we will assume $q_0\in (0,\,1)$.

From~\eqref{eq:Eg_Ef} we see immediately that a condition for aggregation/consensus in the deterministic model~\eqref{eq:Boltzmann_1} is
\begin{equation}
	0<\lambda<\sqrt{\frac{q_0(1-q_0)}{\Var{\theta}}}.
	\label{eq:consensus_g}
\end{equation}
To pursue the analysis of the stochastic model~\eqref{eq:Boltzmann_2} we consider instead two particular but sufficiently representative cases of the probability distribution $\rho$ of $\theta$.

\subsubsection{Uniformly distributed $\theta$}
\label{sect:unif_theta}
If we take $\theta\sim\mathcal{U}(\rangeth)$ with $\rangeth=[-1,\,1]$ and $\Var{\theta}=\frac{1}{3}$ then $\rho(\theta)=\frac{1}{2}\mathbbm{1}_{[-1,\,1]}(\theta)$ and we can compute:
\begin{align*}
	E_{\bar{f}}(t)=\int_{\rangeth}E_f(t;\,\theta)\rho(\theta)\,d\theta
		&= \frac{1}{2}e^{2q_0(q_0-1)t}\int_{-1}^{1}e^{2\lambda\left(\lambda\theta^2+(2q_0-1)\theta\right)t}\,d\theta,
\intertext{which, completing the square at the exponent in the integral, gives}
	&= \frac{1}{4\lambda}\sqrt{\frac{\pi}{2t}}e^{-t/2}\left[\Erfi{\xi_{+}\sqrt{t}}-\Erfi{\xi_{-}\sqrt{t}}\right],
\end{align*}
where $\Erfi{x}:=\frac{2}{\sqrt{\pi}}\int_0^xe^{y^2}\,dy$ is the imaginary error function and where we have denoted
$$ \xi_{\pm}:=\pm\sqrt{2}\lambda+\frac{2q_0-1}{\sqrt{2}}. $$
Since $\Erfi{x}\sim\frac{1}{\sqrt{\pi}x}e^{x^2}$ for $\abs{x}\to\infty$, we deduce
\begin{align*}
	E_{\bar{f}}(t) &\sim \frac{1}{4\sqrt{2}\lambda t}e^{-t/2}\left(\frac{e^{\xi_{+}^2t}}{\xi_{+}}
		-\frac{e^{\xi_{-}^2t}}{\xi_{-}}\right) \\
	&\phantom{\sim}= \pm\frac{1}{4\sqrt{2}\lambda\xi_{\pm}t}e^{\left(-\frac{1}{2}+\xi_{\pm}^2\right)t}
		\left(1-\frac{\xi_{\pm}}{\xi_{\mp}}e^{-\left(\xi_{\pm}^2-\xi_{\mp}^2\right)t}\right) \qquad (t\to+\infty)
\end{align*}
with $\xi_{+}^2-\xi_{-}^2=4\lambda(2q_0-1)$. Writing
\begin{equation}
	E_{\bar{f}}(t)\sim
	\begin{cases}
		-\dfrac{1}{4\sqrt{2}\lambda\xi_{-}t}e^{\left(-\frac{1}{2}+\xi_{-}^2\right)t}
			\left(1-\dfrac{\xi_{-}}{\xi_{+}}e^{4\lambda(2q_0-1)t}\right) & \text{if } 0<q_0\leq\dfrac{1}{2} \\[5mm]
		\dfrac{1}{4\sqrt{2}\lambda\xi_{+}t}e^{\left(-\frac{1}{2}+\xi_{+}^2\right)t}
			\left(1-\dfrac{\xi_{+}}{\xi_{-}}e^{-4\lambda(2q_0-1)t}\right) & \text{if } \dfrac{1}{2}<q_0<1
	\end{cases}
	\label{eq:Ef_bar}
\end{equation}
we see that $E_{\bar{f}}(t)\to 0$ when $t\to+\infty$ provided
\begin{equation*}
	\begin{cases}
		-\dfrac{1}{2}+\xi_{-}^2\leq 0 & \text{if } 0<q_0\leq\dfrac{1}{2} \\[3mm]
		-\dfrac{1}{2}+\xi_{+}^2\leq 0 & \text{if } \dfrac{1}{2}<q_0<1,
	\end{cases}
\end{equation*}
which produces the following condition for aggregation/consensus in the stochastic model~\eqref{eq:Boltzmann_2}:
\begin{equation}
	0<\lambda\leq\min\{q_0,\,1-q_0\}.
	\label{eq:consensus_f}
\end{equation}

\begin{figure}[!t]
\centering
\includegraphics[width=0.6\textwidth]{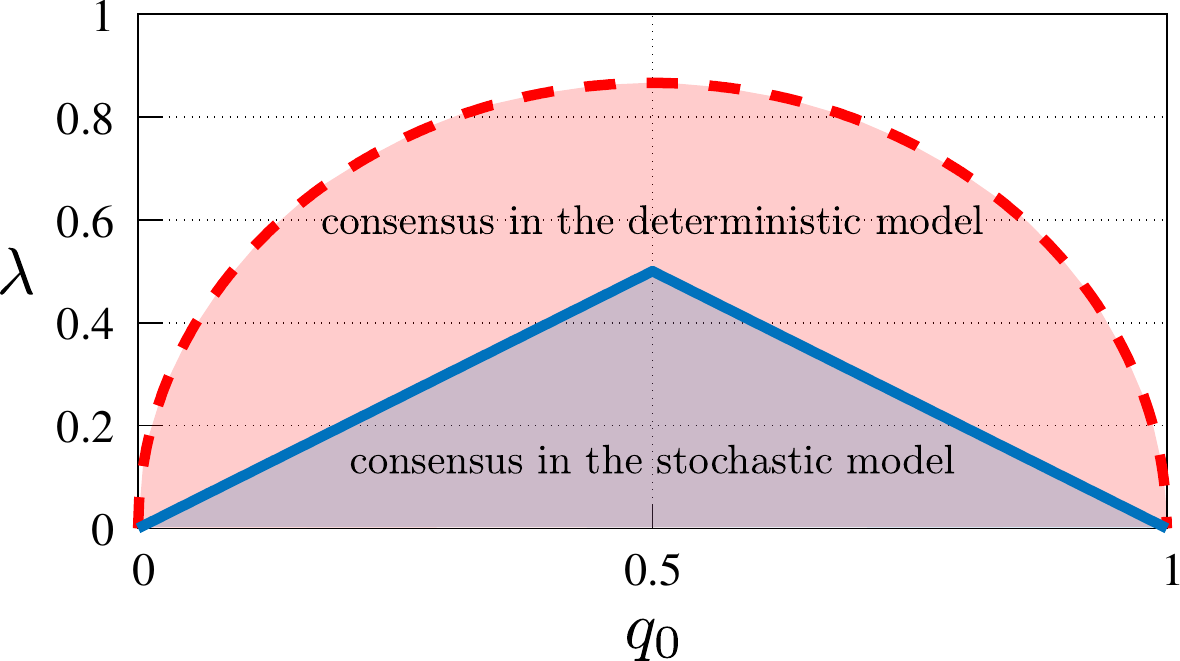}
\caption{Aggregation/consensus regions in the case $\theta\sim\mathcal{U}(-1,\,1)$. Red: condition~\eqref{eq:consensus_g} for the deterministic model with $\Var{\theta}=\frac{1}{3}$. Blue: condition~\eqref{eq:consensus_f} for the stochastic model.}
\label{fig:consensus_cond}
\end{figure}

Figure~\ref{fig:consensus_cond} shows that consensus in the stochastic model happens in a subregion of the plane $(q_0,\,\lambda)$ contained in the consensus region of the deterministic model. This indicates that a stronger restriction is needed on $\lambda$ to reach aggregation/consensus in the former than in the latter.

Furthermore, from~\eqref{eq:Eg_Ef} and~\eqref{eq:Ef_bar} we discover
$$ \frac{E_{\bar{f}}(t)}{E_g(t)}\sim\frac{C(t)}{t}e^{2\left((1-\Var{\theta})\lambda^2+\abs{2q_0-1}\right)t}
	=\frac{C(t)}{t}e^{2\left(\frac{2}{3}\lambda^2+\abs{2q_0-1}\right)t} \qquad (t\to+\infty), $$
where $C(t)$ tends asymptotically to a non-zero finite value. Thus it results invariably $E_g=o(E_{\bar{f}})$ when $t\to+\infty$, meaning that the rate of convergence to aggregation/consensus of the deterministic model is higher than the corresponding rate of the stochastic model.

\begin{remark}
If the space of the microscopic states $V$ is bounded then a further restriction has to be imposed on $\lambda$ in order to guarantee $v^\ast,\,w^\ast\in V$ given $v,\,w\in V$. From~\eqref{eq:binary_symmetric} it can be checked that the bounds on the post-interaction states are never violated if $0\leq q(\theta)\leq 1$, which, taking~\eqref{eq:q} into account together with $\theta\in [-1,\,1]$, implies $0<\lambda\leq\min\{q_0,\,1-q_0\}$, i.e., the same condition as~\eqref{eq:consensus_f}. Therefore we conclude that if $V$ is bounded then the rules~\eqref{eq:binary_symmetric}-\eqref{eq:q} lead always to aggregation/consensus in both the deterministic and the stochastic model, however with a faster convergence in the former case.
\end{remark}

\subsubsection{Normally distributed $\theta$}
If we take $\theta\sim\mathcal{N}(0,\,1)$ with $\rangeth=\R$ then $\rho(\theta)=\frac{1}{\sqrt{2\pi}}e^{-\theta^2/2}$ and we can compute:
\begin{align*}
	E_{\bar{f}}(t)=\int_{\rangeth}E_f(t;\,\theta)\rho(\theta)\,d\theta
		&= \frac{1}{\sqrt{2\pi}}e^{2q_0(q_0-1)t}\int_{-\infty}^{+\infty}e^{2\lambda\left(\lambda\theta^2+(2q_0-1)\theta\right)t}e^{-\theta^2/2}\,d\theta.
\intertext{First, we notice that the integral converges only if $t<\frac{1}{4\lambda^2}$. Under this assumption, completing the square at the exponent in the integral we obtain}
	&= \frac{1}{\sqrt{1-4\lambda^2t}}\exp\left(2q_0\left(q_0-1\right)t+\frac{2\lambda^2\left(2q_0-1\right)^2}{1-4\lambda^2t}t^2\right),
\end{align*}
whence we see that $E_{\bar{f}}\to+\infty$ for $t\uparrow\frac{1}{4\lambda^2}$ for all values of $q_0$ and $\lambda$. Therefore in this case we cannot expect convergence to aggregation/consensus on average in the stochastic model as a side-effect of the non-zero tails of the distribution of $\theta$, cf.~\cite{albi2015MPE}.

\section{Fokker-Planck asymptotics and steady states}
\label{sect:Fokker-Planck}
In order to gain more detailed insights into the large time behaviour of the solutions to~\eqref{eq:Boltzmann_1},~\eqref{eq:Boltzmann_2}, especially when also stochastic fluctuations are included in the microscopic interaction rules, one can resort to particular asymptotics of the Boltzmann equation which provide good approximations of the original steady states. A celebrated one, called the \emph{quasi-invariant interaction limit}~\cite{toscani2006CMS}, results in \emph{Fokker-Planck equations} obtained through scaling techniques which are reminiscent of the classical grazing collision limit~\cite{furioli2017M3AS,pareschi2013BOOK}.

In the context of this paper, the interest is in comparing such asymptotics for the deterministic and the stochastic models. To this purpose, let us consider the following generalisation of the binary interactions~\eqref{eq:binary_general}:
\begin{align}
	\begin{aligned}[c]
		v^\ast &= v+\gamma\left[(p_1(\theta)-1)v+q_1(\theta)w\right]+D(v)\eta \\
		w^\ast &= w+\gamma\left[p_2(\theta)v+(q_2(\theta)-1)w\right]+D(w)\eta,
	\end{aligned}
	\label{eq:binary_general-gamma}
\end{align}
where $\gamma>0$ is a given constant, $\eta$ is a random variable with zero mean, $\Var{\eta}=\sigma^2<+\infty$ and bounded third order moment, $D:V\to\R_+$ is a local diffusion coefficient possibly depending on the microscopic state itself. We further assume that the random variables $\theta$ and $\eta$ are independent. Notice that from~\eqref{eq:binary_general-gamma} we recover~\eqref{eq:binary_general} in the special case $\gamma=1$ and $D=0$.

Let us introduce the time scale $\tau:=\gamma t$ and the corresponding scaled distribution functions:
\begin{equation}
	g_\gamma(\tau,\,v):=g(t,\,v)=g(\tau/\gamma,\,v), \qquad
		f_\gamma(\tau,\,v;\,\theta):=f(t,\,v;\,\theta)=f(\tau/\gamma,\,v;\,\theta)
	\label{eq:scaled_g-f}
\end{equation}
for the deterministic and stochastic models, respectively. From~\eqref{eq:Boltzmann_1},~\eqref{eq:Boltzmann_2} they are readily seen to satisfy
\begin{align}
	\begin{aligned}[b]
		\frac{d}{d\tau} & \int_V\varphi(v)g_\gamma(\tau,\,v)\,dv \\
		&= \frac{1}{2\gamma}\iint_{V^2}\left(\bracket{\ave{\varphi(v^\ast)+\varphi(w^\ast)}}
			-\varphi(v)-\varphi(w)\right)g_\gamma(\tau,\,v)g_\gamma(\tau,\,w)\,dv\,dw
	\end{aligned}
	\label{eq:Boltzmann_1-gamma}
\end{align}
and
\begin{align}
	\begin{aligned}[b]
		\frac{d}{d\tau} & \int_V\varphi(v)f_\gamma(\tau,\,v;\,\theta)\,dv \\
		&= \frac{1}{2\gamma}\iint_{V^2}\left(\bracket{\varphi(v^\ast)+\varphi(w^\ast)}-\varphi(v)-\varphi(w)\right)
			f_\gamma(\tau,\,v;\,\theta)f_\gamma(\tau,\,w;\,\theta)\,dv\,dw,
	\end{aligned}
	\label{eq:Boltzmann_2-gamma}
\end{align}
where $\bracket{\cdot}$ denotes a further average with respect to the distribution of $\eta$.

The aforementioned quasi-invariant interaction limit consists in exploring the regime $\gamma\to 0^+$, $\sigma^2\to 0^+$ along with the assumption that the ratio $\sigma^2/\gamma$ tends simultaneously to a finite nonzero value, say $1$ for simplicity. Then the scaled distribution functions $g_\gamma$, $f_\gamma$ can be shown, cf.~\cite{toscani2006CMS}, to converge, up to subsequences, to distribution functions $g=g(\tau,\,v)$, $f=f(\tau,\,v;\,\theta)$ which are weak solutions to the following Fokker-Planck equations:
\begin{equation}
	\partial_\tau g+\partial_v\left[\left(\int_V\ave{P(v,\,w;\,\theta)}g(\tau,\,w)\,dw\right)g\right]
			=\frac{1}{2}\partial^2_v\left(D^2(v)g\right)
	\label{eq:FP_g}
\end{equation}
and
\begin{equation}
	\partial_\tau f+\partial_v\left[\left(\int_V P(v,\,w;\,\theta)f(\tau,\,w;\,\theta)\,dw\right)f\right]
		=\frac{1}{2}\partial^2_v\left(D^2(v)f\right),
	\label{eq:FP_f}
\end{equation}
where we have set for brevity
\begin{equation}
	P(v,\,w;\,\theta):=\frac{1}{2}\left[(p_1(\theta)+q_2(\theta)-2)v+(p_2(\theta)+q_1(\theta))w\right].
	\label{eq:P}
\end{equation}
Interestingly, owing to~\eqref{eq:scaled_g-f} the large time trends of $g(t,\,v)$ and $f(t,\,v;\,\theta)$ are well approximated by those of $g(\tau,\,v)$ and $f(\tau,\,v;\,\theta)$ satisfying~\eqref{eq:FP_g},~\eqref{eq:FP_f}.

\begin{remark} \label{rem:FP_mean}
Taking $\varphi(v)=v$ in~\eqref{eq:Boltzmann_1-gamma},~\eqref{eq:Boltzmann_2-gamma} reveals that, since $\eta$ has zero mean, the evolution of $m_{g_\gamma}(\tau)$, $m_{f_\gamma}(\tau;\,\theta)$ coincides for every $\gamma>0$, hence also in the limit $\gamma\to 0^+$, with the one predicted by the interaction models~\eqref{eq:Boltzmann_1},~\eqref{eq:Boltzmann_2}. Conversely, the evolution of $E_{g_\gamma}(\tau)$, $E_{f_\gamma}(\tau;\,\theta)$ is in general different if $\gamma\neq 1$, thus in particular in the limit $\gamma\to 0^+$.
\end{remark}

\subsection{The inelastic Kac model}\label{sec:ineKac}
In~\cite{pulvirenti2004JSP} Pulvirenti and Toscani proposed a modification of the Kac model~\eqref{eq:Kac}, then further investigated in~\cite{furioli2012SIMA,furioli2017M3AS}, which has the form~\eqref{eq:binary_general} with
\begin{equation}
	p_1(\theta)=q_2(\theta)=\cos{\theta}\abs{\cos{\theta}}^p, \qquad
		p_2(\theta)=-q_1(\theta)=\sin{\theta}\abs{\sin{\theta}}^p
	\label{eq:inelastic_Kac}
\end{equation}
and $p\geq 0$ is a fixed parameter. For $p=0$ one recovers precisely~\eqref{eq:Kac} while for $p>0$ one obtains an \emph{inelastic} model such that the energy is in general not conserved in a single interaction, in fact
$$ \frac{1}{2}\left[(v^\ast)^2+(w^\ast)^2\right]=
	\frac{1}{2}\left(\abs{\cos{\theta}}^{2(1+p)}+\abs{\sin{\theta}}^{2(1+p)}\right)(v^2+w^2)
		\leq \frac{1}{2}\left(v^2+w^2\right). $$

The inelastic Kac model can be framed in the binary interaction scheme~\eqref{eq:binary_general-gamma} as
\be\label{eq:binary_Kac_ine}
\begin{split}
	v^\ast &= v+\gamma\left[(\cos{\theta}\abs{\cos{\theta}}^p-1)v-\sin{\theta}\abs{\sin{\theta}}^pw\right]+D\eta \\
	w^\ast &= w+\gamma\left[\sin{\theta}\abs{\sin{\theta}}^pv+(\cos{\theta}\abs{\cos{\theta}}^p-1)w\right]+D\eta,
\end{split}
\ee
where, since $V=\R$, we have assumed a constant diffusion coefficient $D>0$. This new model is still dissipative, at least in the absence of stochastic fluctuations ($D=0$), if $\gamma$ is sufficiently small, indeed
\begin{align*}
	\frac{1}{2}\bracket{(v^\ast)^2+(w^\ast)^2} &= \frac{1}{2}\Bigl[\gamma^2\left(\abs{\cos{\theta}}^{2(1+p)}+\abs{\sin{\theta}}^{2(1+p)}\right)+
		2\gamma(1-\gamma)\cos{\theta}\abs{\cos{\theta}}^p \\
	&\phantom{=} +(1-\gamma)^2\Bigr](v^2+w^2)+\sigma^2D^2 \\
	&\leq \frac{1}{2}\left[1-2\gamma(1-\gamma)\left(1-\cos{\theta}\abs{\cos{\theta}}^p\right)\right](v^2+w^2)+\sigma^2 D
\end{align*}
for $\gamma\leq 1$. It is therefore interesting to investigate the asymptotic behaviour of its solution in the quasi-invariant interaction regime.
	
Owing to~\eqref{eq:inelastic_Kac}, the function $P$ featuring in~\eqref{eq:FP_g},~\eqref{eq:FP_f} is $P(v,\,w;\,\theta)=\left(\cos{\theta}\abs{\cos{\theta}}^p-1\right)v$. In particular, recalling that $\theta\sim\mathcal{U}(0,\,2\pi)$,
$$ \ave{P(v,\,w;\,\theta)}=\left(\frac{1}{2\pi}\int_0^{2\pi}\cos{\theta}\abs{\cos{\theta}}^p\,d\theta-1\right)v=-v, $$
therefore the deterministic and stochastic Fokker-Planck models read respectively
$$ \partial_\tau g-\partial_v(vg)=\frac{D^2}{2}\partial^2_v g, \qquad
	\partial_\tau f-(1-\cos{\theta}\abs{\cos{\theta}}^p)\partial_v(vf)=\frac{D^2}{2}\partial^2_vf. $$
They can be solved explicitly at the steady state to obtain the following asymptotic distribution functions $g^\infty=g^\infty(v)$ and $f^\infty=f^\infty(v;\,\theta)$:
$$ g^\infty(v)=\frac{1}{D\sqrt{\pi}}e^{-v^2/D^2}, \qquad
	f^\infty(v;\,\theta)=\frac{\sqrt{1-\cos{\theta}\abs{\cos{\theta}}^p}}{D\sqrt{\pi}}e^{-(1-\cos{\theta}\abs{\cos{\theta}}^p)v^2/D^2}. $$

Consistently with Remark~\ref{rem:FP_mean} and with Section~\ref{sect:Kac}, we observe that $m_{g^\infty}=m_{f^\infty}(\theta)=0$ for all $\theta\in (0,\,2\pi)$, hence also $m_{\bar{f}^\infty}=0$. Concerning the energies we find instead
\be\label{eq:E_K_ine_1}
	E_{g^\infty}=\frac{D^2}{2}, \qquad E_{f^\infty}(\theta)=\frac{D^2}{2(1-\cos{\theta}\abs{\cos{\theta}}^p)}, 
\ee
whence
\be\label{eq:E_K_ine_2}
	E_{\bar{f}^\infty}=\ave{E_{f^\infty}(\theta)}=\frac{D^2}{4\pi}\int_0^{2\pi}\frac{1}{1-\cos{\theta}\abs{\cos{\theta}}^p}\,d\theta=+\infty 
\ee
because $1-\cos{\theta}\abs{\cos{\theta}}^p\sim\frac{1+p}{2}\theta^2$ for $\theta\to 0^+$ and $1-\cos{\theta}\abs{\cos{\theta}}^p\sim\frac{1+p}{2}(\theta-2\pi)^2$ for $\theta\to 2\pi^-$.

The finite energy predicted asymptotically by the deterministic model is a consequence of the balance between the dissipation and the diffusion taking place on average in each binary interaction as observed above. Remarkably, the stochastic model predicts instead asymptotically an infinite energy on average, due to that for $\theta$ close to $0,\,2\pi$ the dynamics~\eqref{eq:binary_general-gamma} are essentially dominated by the constant diffusion.

\subsection{Aggregation/consensus models with diffusion}
If we consider the binary interaction scheme~\eqref{eq:binary_general-gamma} in the setting introduced in Section~\ref{sect:consensus} we find
\begin{equation}
	v^\ast=v+\gamma q(\theta)(w-v)+D(v)\eta, \qquad
		w^\ast=w+\gamma q(\theta)(v-w)+D(w)\eta
	\label{eq:binary_symmetric_diff}
\end{equation}
with $q(\theta)=q_0+\lambda\theta$. From Remark~\ref{rem:FP_mean} we know that $m_g(\tau)$ and $m_f(\tau;\,\theta)$ evolve like in Section~\ref{sect:consensus}, i.e. $m_g(\tau)=m_f(\tau;\,\theta)=0$ for all $\tau>0$ and all $\theta\in\rangeth$ if they vanish for $\tau=0$. From~\eqref{eq:Boltzmann_1-gamma},~\eqref{eq:Boltzmann_2-gamma} with $\varphi(v)=v^2$ we discover that for $\gamma\to 0^+$ it results
\begin{align*}
	& \frac{d}{d\tau}E_{g}(\tau)=-2q_0E_{g}(\tau)+\int_VD^2(v)g(\tau,\,v)\,dv \\
	& \frac{d}{d\tau}E_{f}(\tau;\,\theta)=-2q(\theta)E_{f}(\tau;\,\theta)+\int_VD^2(v)f(\tau,\,v;\,\theta)\,dv,
\end{align*}
whence, averaging the second equation with respect to $\theta$,
$$ \frac{d}{d\tau}E_{\bar{f}}(\tau)=-2q_0E_{\bar{f}}(\tau)
	-2\lambda\ave{\theta E_{f}(\tau;\,\theta)}+\int_VD^2(v)\bar{f}(\tau,\,v)\,dv. $$
Therefore, while in the deterministic model there is an asymptotic balance between dissipative and diffusive effects, in the stochastic model such a trend is perturbed by an additional contribution depending on the variability of $\theta$.

More in general, the deterministic and stochastic Fokker-Planck equations~\eqref{eq:FP_g},~\eqref{eq:FP_f} take now the forms
\begin{equation}
	\partial_\tau g-q_0\partial_v(vg)=\frac{1}{2}\partial^2_v\left(D^2(v)g\right), \qquad
		\partial_\tau f-q(\theta)\partial_v(vf)=\frac{1}{2}\partial^2_v\left(D^2(v)f\right),
	\label{eq:FP_consensus}
\end{equation}
respectively.

\subsubsection{Models on the real line with constant diffusion}
\label{sect:const_diff_R}
Let us consider at first $V=\R$ with $D>0$ constant. The asymptotic solutions to~\eqref{eq:FP_consensus} write
$$ g^\infty(v)=\frac{\sqrt{q_0}}{D\sqrt{\pi}}e^{-q_0v^2/D^2}, \qquad
	f^\infty(v;\,\theta)=\frac{\sqrt{q_0+\lambda\theta}}{D\sqrt{\pi}}e^{-(q_0+\lambda\theta)v^2/D^2} $$
provided $q_0+\lambda\theta>0$, which can be achieved if $\rangeth$ is bounded from the left. In particular, the asymptotic energies are
$$ E_{g^\infty}=\frac{D^2}{2q_0}, \qquad E_{f^\infty}(\theta)=\frac{D^2}{2(q_0+\lambda\theta)}. $$

Considering $\theta\sim\mathcal{U}(-1,\,1)$ with $0<\lambda<q_0$, we obtain the following explicit expression of the average asymptotic distribution function of the stochastic model:
$$ \bar{f}^\infty(v)=\frac{1}{2\pi}\int_{-1}^{1}f^\infty(v;\,\theta)\,d\theta=
	\frac{1}{v^2}\left(C_1e^{-(q_0-\lambda)v^2/D^2}+C_2e^{-(q_0+\lambda)v^2/D^2}+\frac{C_3(v)}{v}\right) $$
where
\begin{gather*}
	C_1:=\frac{D}{2\lambda\pi\sqrt{\pi}}\sqrt{q_0-\lambda}, \qquad
		C_2:=-\frac{D}{2\lambda\pi\sqrt{\pi}}\sqrt{q_0+\lambda} \\
	C_3(v):=\frac{D^2}{4\lambda\pi}\left[\Erf{\frac{\sqrt{q_0+\lambda}}{D}v}-\Erf{\frac{\sqrt{q_0-\lambda}}{D}v}\right].
\end{gather*}
Moreover
$$ E_{\bar{f}^\infty}=\ave{E_{f^\infty}(\theta)}
	=\frac{D^2}{4\lambda}\log{\left(\frac{q_0+\lambda}{q_0-\lambda}\right)}
		=\frac{q_0}{2\lambda}\log\left(\frac{1+\frac{\lambda}{q_0}}{1-\frac{\lambda}{q_0}}\right)E_{g^\infty}, $$
which implies
$$ \frac{E_{\bar{f}^\infty}}{E_{g^\infty}}=\frac{q_0}{2\lambda}\log\left(\frac{1+\frac{\lambda}{q_0}}{1-\frac{\lambda}{q_0}}\right)\geq 1 $$
for $0<\frac{\lambda}{q_0}<1$, the equality holding true in the limit $\lambda\to 0^+$. Notice that $E_{\bar{f}^\infty}\to +\infty$ when $\lambda\to q_0^-$.

With respect to the case $\gamma=1$ (without diffusion) discussed in Section~\ref{sect:unif_theta}, we observe that in the quasi invariant interaction regime the stochastic model has asymptotically finite energy on average under the restriction $0<\lambda<q_0$, which for $q_0>\frac{1}{2}$ is weaker than~\eqref{eq:consensus_f}.

\subsubsection{Models on a bounded interval with non-linear diffusion}
\label{sect:non-linear_diff}
We consider now the case of a bounded $V$, say $V=[-1,\,1]$ to fix the ideas, which requires a non-linear diffusion coefficient $D=D(v)$ vanishing at $v=\pm 1$ in order for~\eqref{eq:binary_symmetric_diff} to preserve the bounds $-1\leq v^\ast,\,w^\ast\leq 1$. The asymptotic solutions to~\eqref{eq:FP_consensus} read
\begin{align*}
	& g^\infty(v)=C_1\exp{\left(-2\int\frac{q_0v+D(v)D'(v)}{D^2(v)}\,dv\right)} \\[2mm]
	& f^\infty(v;\,\theta)=C_2(\theta)\exp{\left(-2\int\frac{q(\theta)v+D(v)D'(v)}{D^2(v)}\,dv\right)},
\end{align*}
$C_1,\,C_2(\theta)>0$ being normalisation constants such that $g^\infty$ and $f^\infty(\cdot;\,\theta)$ have unit integral in $V$. They can be further manipulated for particular choices of the function $D$, cf.~\cite{toscani2006CMS}: for instance, if $D(v)\propto 1-v^2$,
\begin{equation}
	g^\infty(v)=\frac{C_1}{\left(1-v^2\right)^2}e^{-q_0/\left(1-v^2\right)}, \qquad
		f^\infty(v;\,\theta)=\frac{C_2(\theta)}{\left(1-v^2\right)^2}e^{-q(\theta)/\left(1-v^2\right)}.
	\label{eq:gf_non-lin_diff}
\end{equation}
Nevertheless, the explicit determination of the constants $C_1$, $C_2(\theta)$ as well as of the moments of $g^\infty$, $f^\infty(\cdot;\,\theta)$ with respect to $v$ is in general not possible by standard analytical methods. Therefore in the next section we will push ahead with the study of these cases by means of accurate numerical methods.

\section{Numerical examples}
\label{sect:numerics}
In this section we investigate numerically the deterministic and stochastic trends of the Boltz\-mann and Fokker-Planck equations~\eqref{eq:Boltzmann_1}-\eqref{eq:Boltzmann_2} and~\eqref{eq:FP_g}-\eqref{eq:FP_f}, respectively. Some references on popular numerical methods, such as stochastic collocation, stochastic Galerkin schemes, multi-level Monte Carlo schemes, are e.g.~\cite{hu2016JCP,hu2017CHAPTER,mishra2013CHAPTER,xiu2010BOOK,xiu2005SISC,zhu2017PREPRINT}.

Among the numerical methods for quantifying the uncertainty in PDEs like~\eqref{eq:Boltzmann_2} or~\eqref{eq:FP_f}, here we consider collocation methods. Given a discretisation $\{\theta_k\}_{k=0}^M$ of the set $\rangeth$, they consist in solving $M+1$ equations of the form~\eqref{eq:Boltzmann_2} or~\eqref{eq:FP_f} with $\theta=\theta_k$, for $k=0,\,\dots,\,M$. As a result, one gets an ensemble of $M+1$ solutions $\{f(\cdot,\,\cdot;\,\theta_k)\}_{k=0}^{M}$, which can be post-processed in order to obtain statistical information such as the average distribution function $\bar{f}$ and its variance with respect to $\theta$. The collocation nodes $\theta_k$ are typically chosen according to Gaussian quadrature rules, consistently with the probability distribution of $\theta$. The advantage of this approach is that all the properties of the deterministic schemes used to solve the $M+1$ equations with $\theta=\theta_k$ are automatically preserved.

\subsection{Boltzmann-type models}
We start by dealing with the deterministic and stochastic Boltzmann-type equations~\eqref{eq:Boltzmann_1},~\eqref{eq:Boltzmann_2}, which we solve by means of direct simulation Monte Carlo methods (for each collocation node $\theta=\theta_k$ in the case of~\eqref{eq:Boltzmann_2}). In particular, we apply the Nanbu-Babovsky scheme, see~\cite{pareschi2001ESAIMP,pareschi2013BOOK} for an introduction. In each numerical test we consider $10^5$ particles.

\begin{figure}[!t]
\centering
\includegraphics[scale=0.5]{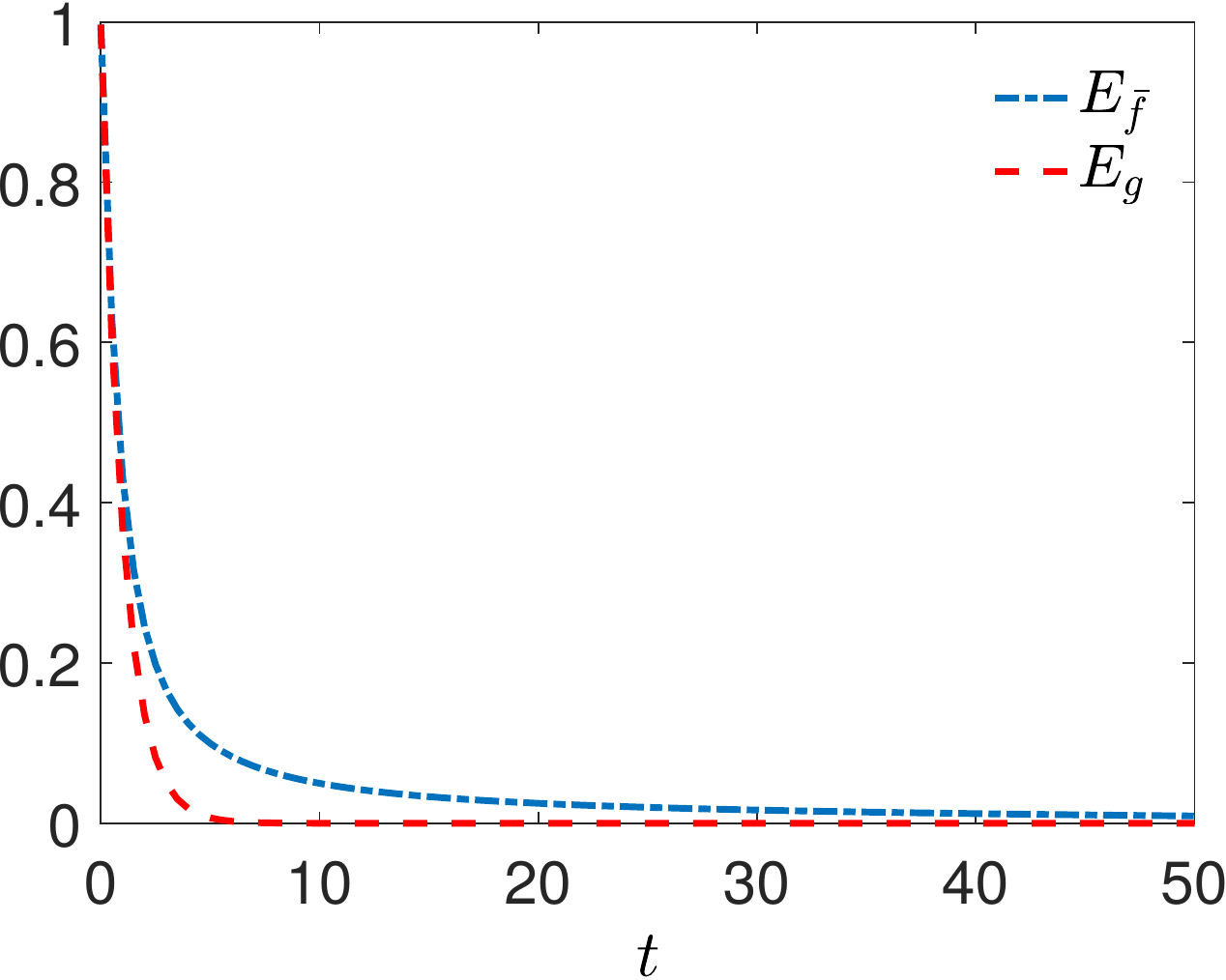}
\caption{\emph{Example 1a}. Decay of the energies produced by the deterministic and stochastic Boltzmann-type models~\eqref{eq:Boltzmann_1},~\eqref{eq:Boltzmann_2} in the time interval $[0,\,50]$ under the binary interaction rules~\eqref{eq:binary_symmetric}-\eqref{eq:q} with $q_0=\lambda=0.5$ and $\theta\sim\mathcal{U}(-1,\,1)$. For the direct Monte Carlo approximation of~\eqref{eq:Boltzmann_2} we used a collocation method with $M=10$ collocation nodes.}
\label{fig:Boltzmann_energies}
\end{figure}

\begin{figure}[!t]
\centering
\subfigure[$\bar{f}(t,\,v)$]{\includegraphics[scale=0.5]{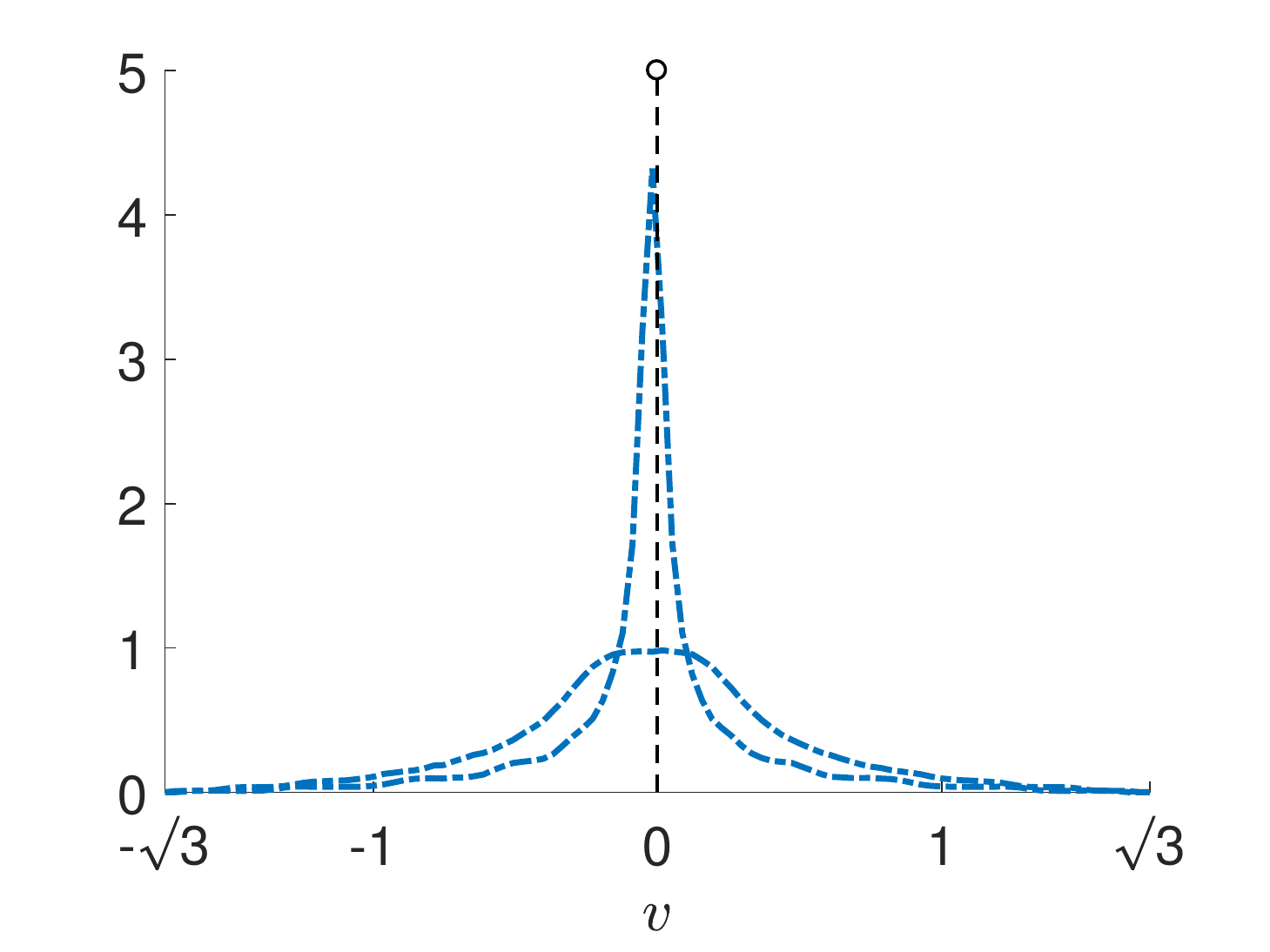}}
\subfigure[$g(t,\,v)$]{\includegraphics[scale=0.5]{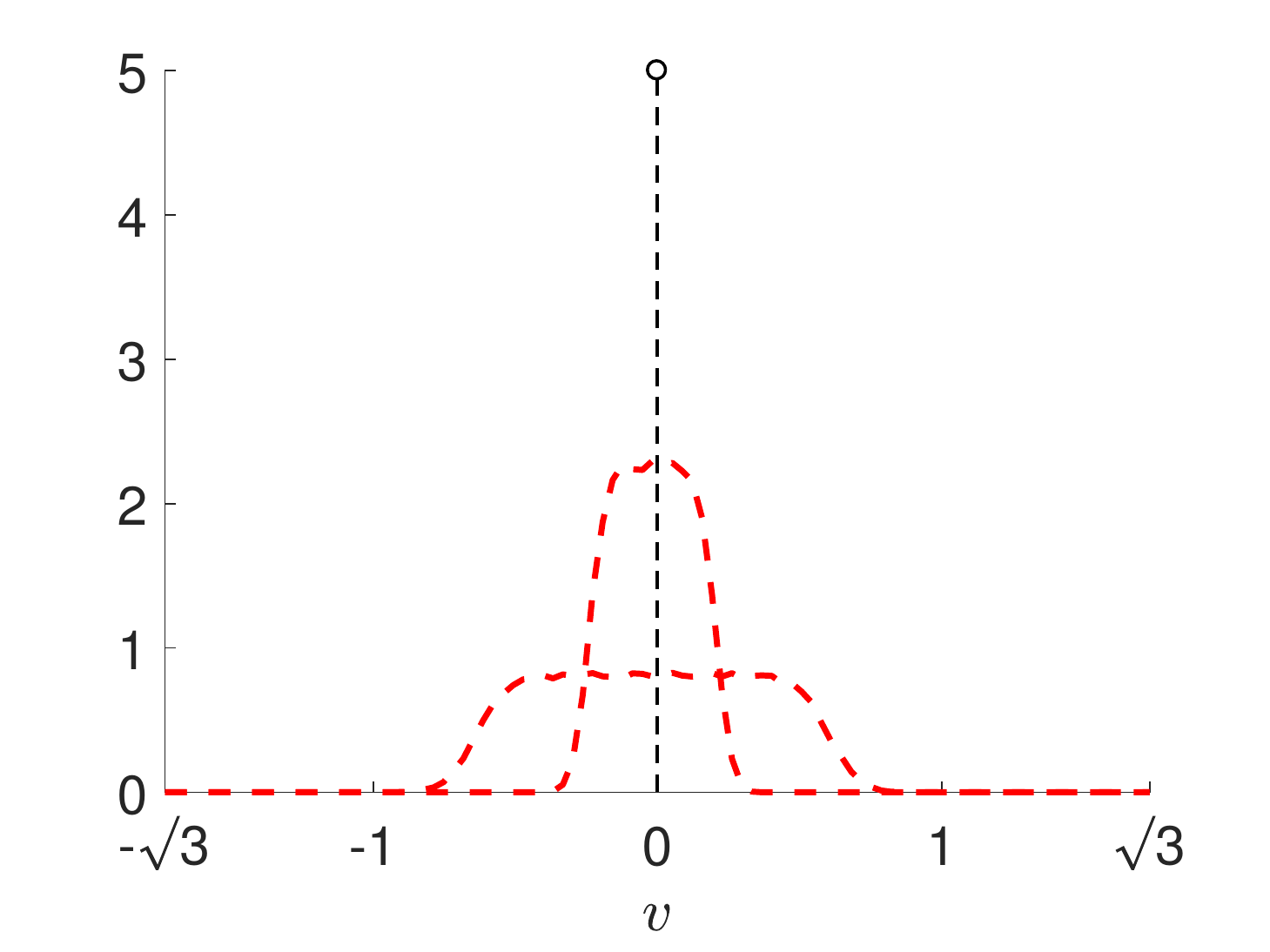}}
\caption{\emph{Example 1a}. Approximation of transient distributions of~\eqref{eq:Boltzmann_1},~\eqref{eq:Boltzmann_2} at the time steps $t=2$ and $t=4$ with binary interaction rules like in Figure~\ref{fig:Boltzmann_energies}. The black dashed vertical line represents the asymptotic Dirac delta centred in the (conserved) null mean.}
\label{fig:Boltzmann_transient}
\end{figure}

\subsubsection{Example 1a: Aggregation/consensus}
To begin with, we focus on the aggregation/consensus model with binary interactions~\eqref{eq:binary_symmetric} with uncertain parameter $q(\theta)$ given by~\eqref{eq:q} and $\theta\sim\mathcal{U}(-1,\,1)$. As far as the numerical approximation of~\eqref{eq:Boltzmann_2} is concerned, we consider a collocation method with $M=10$ Legendre nodes. We recall that, as shown in Figure \ref{fig:consensus_cond}, under several combinations of the parameters $q_0$, $\lambda$ the deterministic and stochastic Boltzmann-type models~\eqref{eq:Boltzmann_1}-\eqref{eq:Boltzmann_2} converge asymptotically to the same stationary state, i.e. the Dirac delta distribution centred in the (conserved) mean of the system.

In Figure~\ref{fig:Boltzmann_energies} we show the evolution in the time interval $[0,\,50]$ of the energies $E_g$ and $E_{\bar{f}}$ for $q_0=\lambda=0.5$. Although with this choice of the parameters both energies are expected to tend asymptotically to zero, cf.~\eqref{eq:consensus_f} and Figure~\ref{fig:consensus_cond}, we observe that the rate of convergence in the deterministic model is higher than the corresponding rate in the stochastic model, consistently with the theoretical findings of Section~\ref{sect:unif_theta}. This behaviour is also confirmed by the transient distribution functions $g(t,\,v)$ and $\bar{f}(t,\,v)$, which are shown in Figure~\ref{fig:Boltzmann_transient} at two successive times $t=2,\,4$.

\begin{figure}[!t]
\centering
\subfigure[]{\includegraphics[scale=0.5]{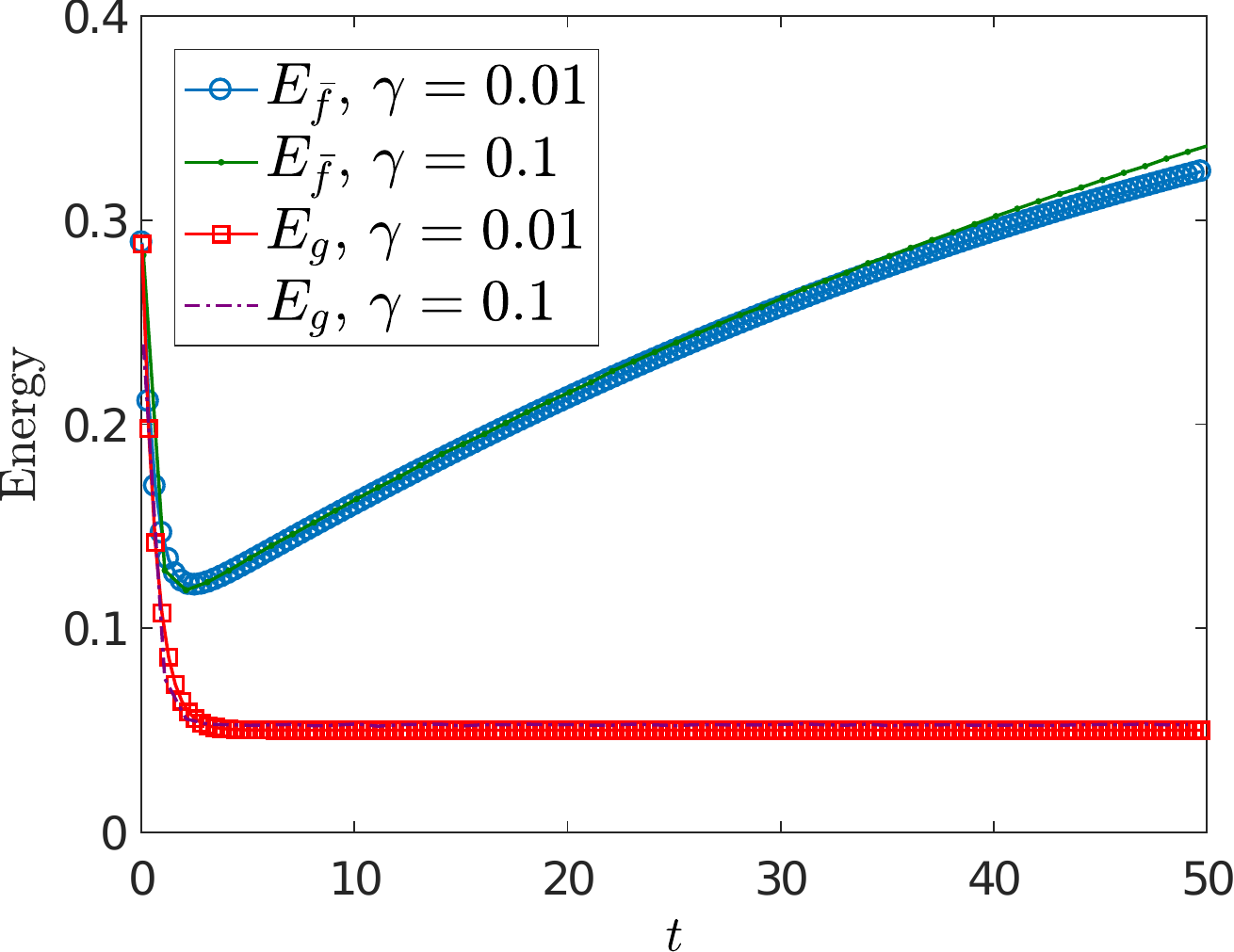}}
\subfigure[]{\includegraphics[scale=0.5]{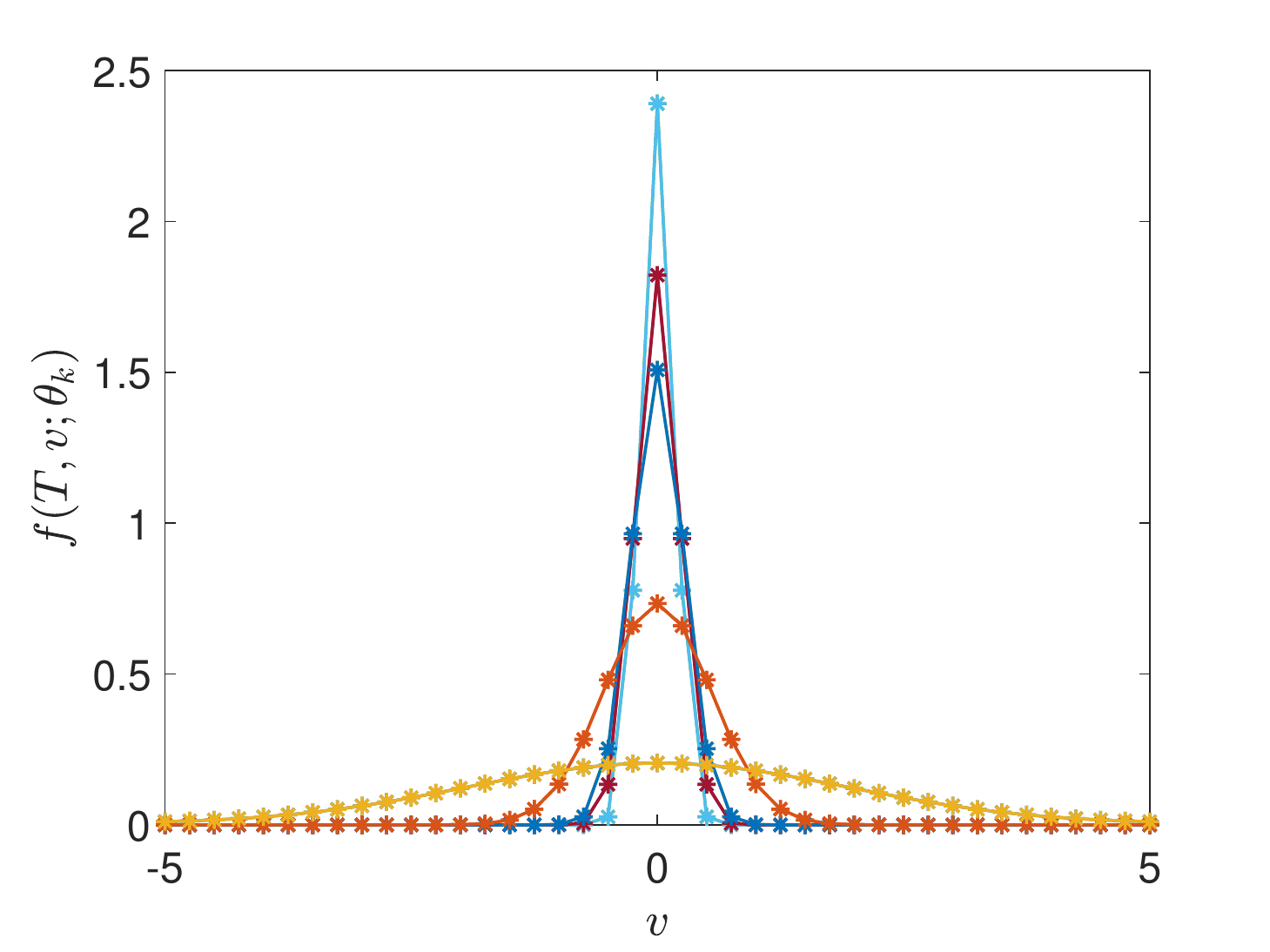}}
\caption{(a) Evolution in the time interval $[0,\,50]$ of the energies of the deterministic and stochastic inelastic Kac model with $p=1$ computed from the scaled Boltzmann equations~\eqref{eq:Boltzmann_1-gamma},~\eqref{eq:Boltzmann_2-gamma} with $\Var{\eta}=\gamma$ and the two values $\gamma=10^{-1}$, $\gamma=10^{-2}$. In addition to that we set $D^2=0.1$ and $\Delta{t}=\gamma$ in the direct Monte Carlo method. (b) Stationary distributions of~\eqref{eq:Boltzmann_2-gamma} computed for different random values of $\theta=\theta_k$, $k=1,\,\dots,\,10$.}
\label{fig:Kac_model}
\end{figure}

\subsubsection{Example 2a: Inelastic Kac model}
We consider now the inelastic Kac model introduced in Section~\ref{sec:ineKac}, defined by the binary interactions~\eqref{eq:binary_Kac_ine} with furthermore $p=1$ and $D^2=0.1$ constant. We solve numerically, by means of a collocation scheme in the frame of direct Monte Carlo methods, the Boltzmann equations~\eqref{eq:Boltzmann_1-gamma},~\eqref{eq:Boltzmann_2-gamma} for small values of $\gamma>0$ and $\Var{\eta}=\gamma$, so as to mimic the limiting behaviour modelled by the Fokker-Planck equations~\eqref{eq:FP_g},~\eqref{eq:FP_f} in the quasi-invariant interaction limit.

In Figure~\ref{fig:Kac_model}a we show the evolution in the time interval $[0,\,50]$ of the energies $E_g$, $E_{\bar{f}}$ for $\gamma=10^{-1}$ and $\gamma=10^{-2}$. We observe that for $\gamma\ll 1$ the asymptotic trends computed analytically from the Fokker-Planck equations, cf.~\eqref{eq:E_K_ine_1},~\eqref{eq:E_K_ine_2}, are indeed captured by the scaled Boltzmann equations. In particular, we notice that $E_g$ converges to the value $D^2/2=0.05$ while $E_{\bar{f}}$ diverges.

In Figure~\ref{fig:Kac_model}b we show, for the stochastic equation~\eqref{eq:Boltzmann_2-gamma}, various approximations of the stationary distribution functions $\{f^\infty(\cdot;\,\theta)\}$ computed at time $\tau=50$ for $M=10$ random collocation nodes $\theta=\theta_k$. We notice that the closer $\theta_k$ to $0$ or $2\pi$ the more smeared the profile of the corresponding $f^\infty(v;\,\theta_k)$, which, consistently with the findings of Section~\ref{sec:ineKac}, is responsible for the blow-up of $E_{\bar{f}^\infty}$ already observed.

\subsection{Fokker-Planck models}
We now consider the deterministic and stochastic models obtained by means of the scaling techniques discussed in Section~\ref{sect:Fokker-Planck}. In particular, we solve numerically the resulting Fokker-Planck equations~\eqref{eq:FP_g},~\eqref{eq:FP_f} taking advantage of \emph{Structure Preserving} (SP) collocation methods, see~\cite{dimarco2017CHAPTER,pareschi2017JSC}, which preserve important structural properties of the exact solutions such as the non-negativity, the entropy dissipation and the large time behaviour. Furthermore they are second order accurate in the transient regime, do not impose any restrictions on the $v$-mesh size and can capture the steady states with an arbitrary accuracy.

In order to be self-consistent, we summarize in the following some features of SP methods in the stochastic collocation setting. To begin with, we observe that the stochastic Fokker-Planck equation~\eqref{eq:FP_f} may be rewritten for all $k=0,\,\dots,\,M$ in flux form as
$$ \partial_\tau f(\tau,\,v;\,\theta_k)=\partial_v\F[f](\tau,\,v;\,\theta_k), $$
where
\be \label{eq:flux_an}
\F[f](\tau,\,v;\,\theta_k):=\C[f](\tau,\,v;\,\theta_k)f(\tau,\,v;\,\theta_k)+\frac{1}{2}D^2(v)\partial_v f(\tau,\,v;\,\theta_k)
\ee
is the flux and
$$ \C[f](\tau,\,v;\,\theta_k):=-\int_V P(v,\,w;\,\theta_k)f(\tau,\,w;\,\theta_k)\,dw+\frac{1}{2}\partial_vD^2(v). $$

Next we introduce a uniform grid $\{v_i\}_{i=1}^{N}\subset V$ with $\Delta{v}:=v_{i+1}-v_i>0$ constant, we denote $v_{i\pm 1/2}:=v_i\pm\frac{1}{2}\Delta{v}$ and we consider the conservative discretisation
$$ \frac{d}{d\tau}f_i^k(\tau)=\frac{\F_{i+1/2}^k[f](\tau)-\F_{i-1/2}^k[f](\tau)}{\Delta{v}}, \qquad
	i=1,\,\dots,\,N, $$
where $f_i^k(\tau)\approx\frac{1}{\Delta{v}}\int_{v_{i-1/2}}^{v_{i+1/2}}f(\tau,\,v;\,\theta_k)\,dv$ denotes the numerical approximation of the cell average. In particular, we choose a numerical flux function of the form
\be\label{eq:F_num}
\F_{i+1/2}^k[f]:=\tilde{\C}^k_{i+1/2}\tilde{f}_{i+1/2}^k+\frac{1}{2}D^2_{i+1/2}\frac{f_{i+1}^k-f_i^k}{\Delta{v}},
\ee
and at the cell interface $i+\frac{1}{2}$ we define
\be\label{eq:F_delta}
\tilde{f}_{i+1/2}^k:=(1-\delta_{i+1/2}^k)f_{i+1}^k+\delta^k_{i+1/2}f_i^k,
\ee
which is a convex linear combination of the values of $f^k$ in the two adjacent cells $i$, $i+1$. Setting in particular
\be
\label{eq:Ctilde_def}
\tilde{\C}_{i+1/2}^k:=\dfrac{D^2_{i+1/2}}{2\Delta v}\int_{v_i}^{v_{i+1}}
	\dfrac{-\int_V P(v,\,w;\,\theta_k)f(\tau,\,w;\,\theta_k)\,dw+\frac{1}{2}\partial_vD^2(v)}{\frac{1}{2}D^2(v)}\,dv
\ee
we obtain explicitly
\be
\label{eq:delta_def}
\delta^k_{i+1/2}=\frac{1}{\lambda_{i+1/2}^k}+\frac{1}{1-\exp(\lambda_{i+1/2}^k)} \quad \text{where} \quad
	\lambda_{i+1/2}^k:=\frac{\Delta{v}\,\tilde{\C}^k_{i+1/2}}{\frac{1}{2}D^2_{i+1/2}}
\ee
and the following result holds, cf.~\cite{pareschi2017JSC}:
\begin{proposition}
The numerical flux function~\eqref{eq:F_num}-\eqref{eq:F_delta} with $\tilde{\C}_{i+1/2}^k$, $\delta_{i+1/2}^k$ defined in~\eqref{eq:Ctilde_def}-\eqref{eq:delta_def} vanishes when the analytical flux~\eqref{eq:flux_an} is equal to zero in the cell $[v_i,\,v_{i+1}]$. Moreover, $\delta_{i+1/2}^k\in [0,\,1]$ for all $i$ and all $k=0,\,\dots,\,M$.
\end{proposition}

This SP scheme offers the following advantages:
\begin{itemize}
\item For linear problems, it preserves the steady state exactly by choosing
$$ \delta_{i+1/2}^{\infty,k}:=\frac{1}{\log(f_i^{\infty,k})}
	+\frac{f_{i+1}^{\infty,k}}{f_{i+1}^{\infty,k}-f_i^{\infty,k}}, \qquad
	\lambda_{i+1/2}^{\infty,k}=\log\left(\dfrac{f_i^{\infty,k}}{f^{\infty,k}_{i+1}}\right). $$
\item Numerical solutions reproduce, for large times, the statistical properties of the exact steady state with high accuracy if suitable high-order quadrature formulas are used to compute~\eqref{eq:Ctilde_def} for all $k=0,\,\dots,\,M$. To indicate the order of accuracy of the SP scheme we use the notation $\text{SP}_j$, where $j>2$ is the order of the employed quadrature method. By $j=\text{G},\,\text{E}$ we mean a Gaussian-type quadrature and the exact integration, respectively. In the following numerical tests we perform Gaussian integrations with $6$ quadrature points in each cell $[v_i,\,v_{i+1}]$.
\item Non-negativity of the numerical solution, without any restrictions on $\Delta{v}$, may be proved for general strong stability preserving and high-order semi-implicit methods up to specific restrictions on the time step $\Delta{\tau}$, see~\cite{dimarco2017CHAPTER,pareschi2017JSC}.
\item The dissipation of the numerical entropy has been proved for specific problems with gradient flow structure~\cite{pareschi2017JSC}.
\end{itemize}

An analogous scheme may be considered for the approximation of $g(\tau,\,v)$ in~\eqref{eq:FP_g}. In this case, the evolution is totally independent of the random input, thus no $\theta$-post-processing is needed.

\medskip

\begin{figure}[!t]
\centering
\subfigure[]{\includegraphics[scale=0.5]{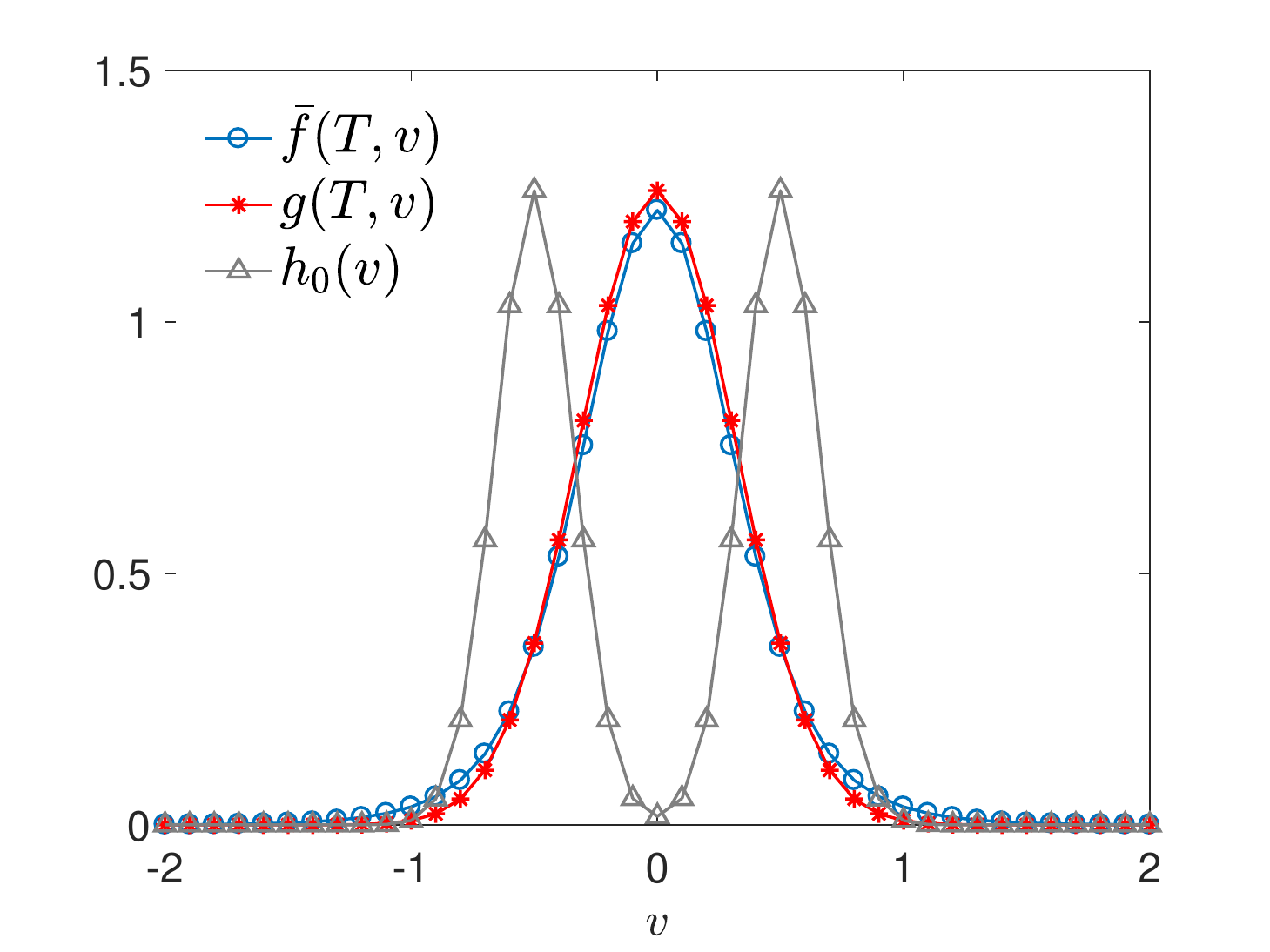}}
\subfigure[]{\includegraphics[scale=0.5]{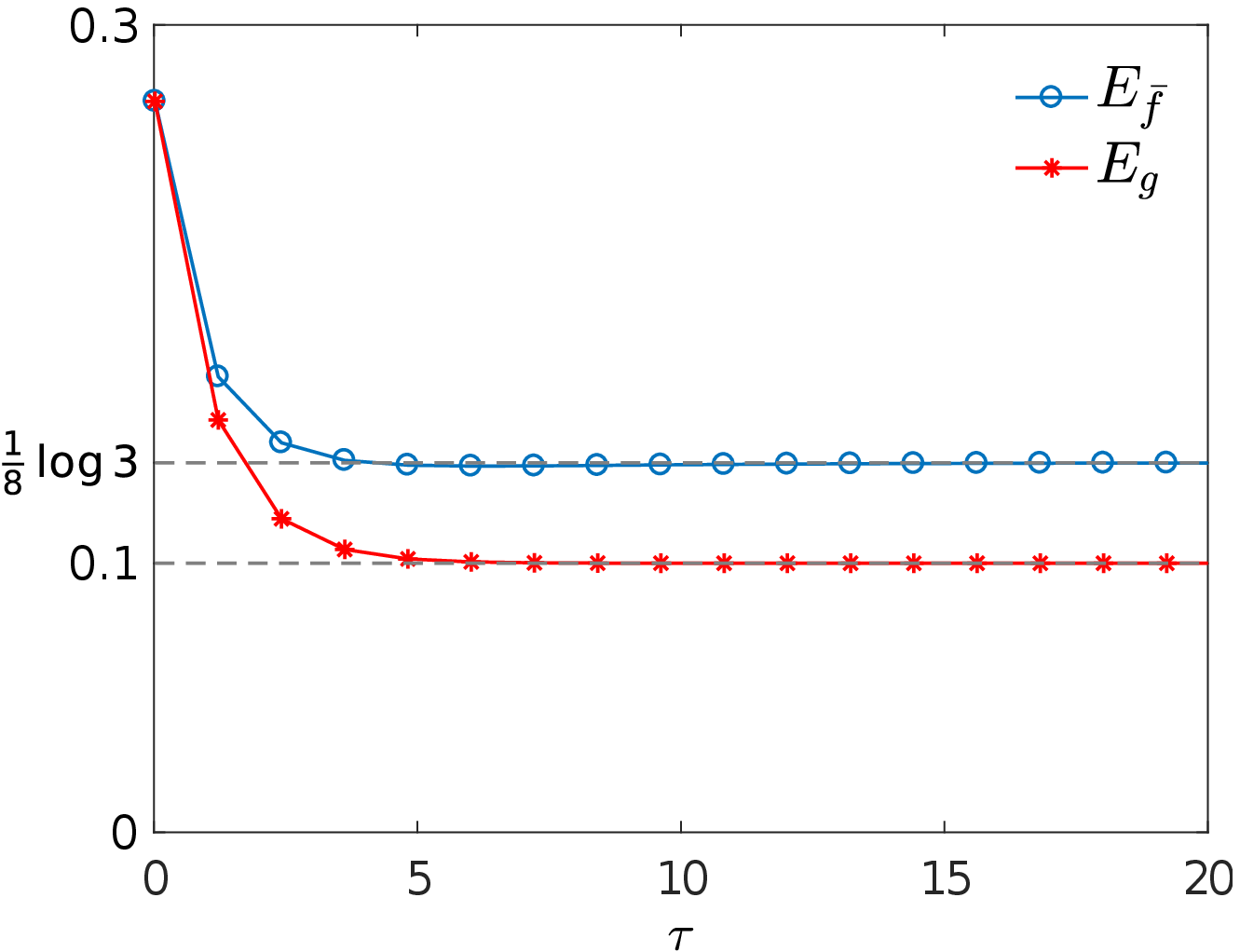}}
\caption{\textit{Example 1b}. (a) Initial (grey) and stationary (red, blue) distributions of the deterministic and the stochastic Fokker-Planck equations in the case of constant diffusion. (b) Evolution in the time interval $[0,\,20]$ of the energies of the two models. For both models we considered the initial distribution $h_0$ in~\eqref{eq:h0}, the bounded computational domain $[-5,\,5]$, $N=101$ grid points and the time step $\Delta{\tau}=2.5\cdot 10^{-4}$. The uncertain interaction parameter~\eqref{eq:q} is such that $q_0=0.5$ and $\lambda=0.4$ while the diffusion coefficient is $D^2=0.1$. The evolutions of $\{f(\tau,\,v;\,\theta_k)\}_{k=0}^{10}$ and of $g(\tau,\,v)$ have been computed with $\text{SP}_\text{E}$ stochastic collocation methods.}
\label{fig:FP_SS}
\end{figure}

In the following tests we invariably consider, for both models~\eqref{eq:FP_g},~\eqref{eq:FP_f}, a deterministic bimodal initial distribution of the form 
\begin{equation}
	g(0,\,v)=f(0,\,v;\,\theta)=h_0(v):=C\left(e^{-20(v-\frac{1}{2})^2}+e^{-20(v+\frac{1}{2})^2}\right),
	\label{eq:h0}
\end{equation}
see Figure~\ref{fig:FP_SS}a (grey line), where $C>0$ is a normalisation constant such that $\int_V h_0(v)\,dv=1$.

\subsubsection{Example 1b: Consensus with constant diffusion}
Let us consider at first~\eqref{eq:FP_g},~\eqref{eq:FP_f} for $$P(v,\,w;\,\theta)=q(\theta)(w-v),\qquad q(\theta)=q_0+\lambda\theta$$ with $\theta\sim\mathcal{U}(-1,\,1)$, and constant diffusion coefficient $D^2=0.1$. This corresponds to the cases analysed in Section~\ref{sect:const_diff_R}. We solve numerically the equations for $v\in [-L,\,L]$, $L=5$, with $N=101$ grid points, hence $\Delta{v}=L/(N-1)=\frac{1}{20}$, while partitioning the time interval $[0,\,T]$, $T=20$, with time step $\Delta{\tau}=\Delta{v}^2/(2L)=2.5\cdot 10^{-4}$.

In Figure~\ref{fig:FP_SS}a we present the numerical solutions $g(T,\,v)$, $\bar{f}(T,\,v)$, which for the chosen value of $T$ can be taken as good approximations of the asymptotic distributions $g^\infty$, $\bar{f}^\infty$, while in Figure~\ref{fig:FP_SS}b we show the time trend of the respective energies. The results are consistent with the analytical steady states described in Section~\ref{sect:const_diff_R}. We stress in particular that, thanks to the constant diffusion and to the properties of the numerical scheme, the large time solutions are described exactly by the numerical approach. 

\begin{figure}[!t]
\centering
\subfigure[]{\includegraphics[scale=0.5]{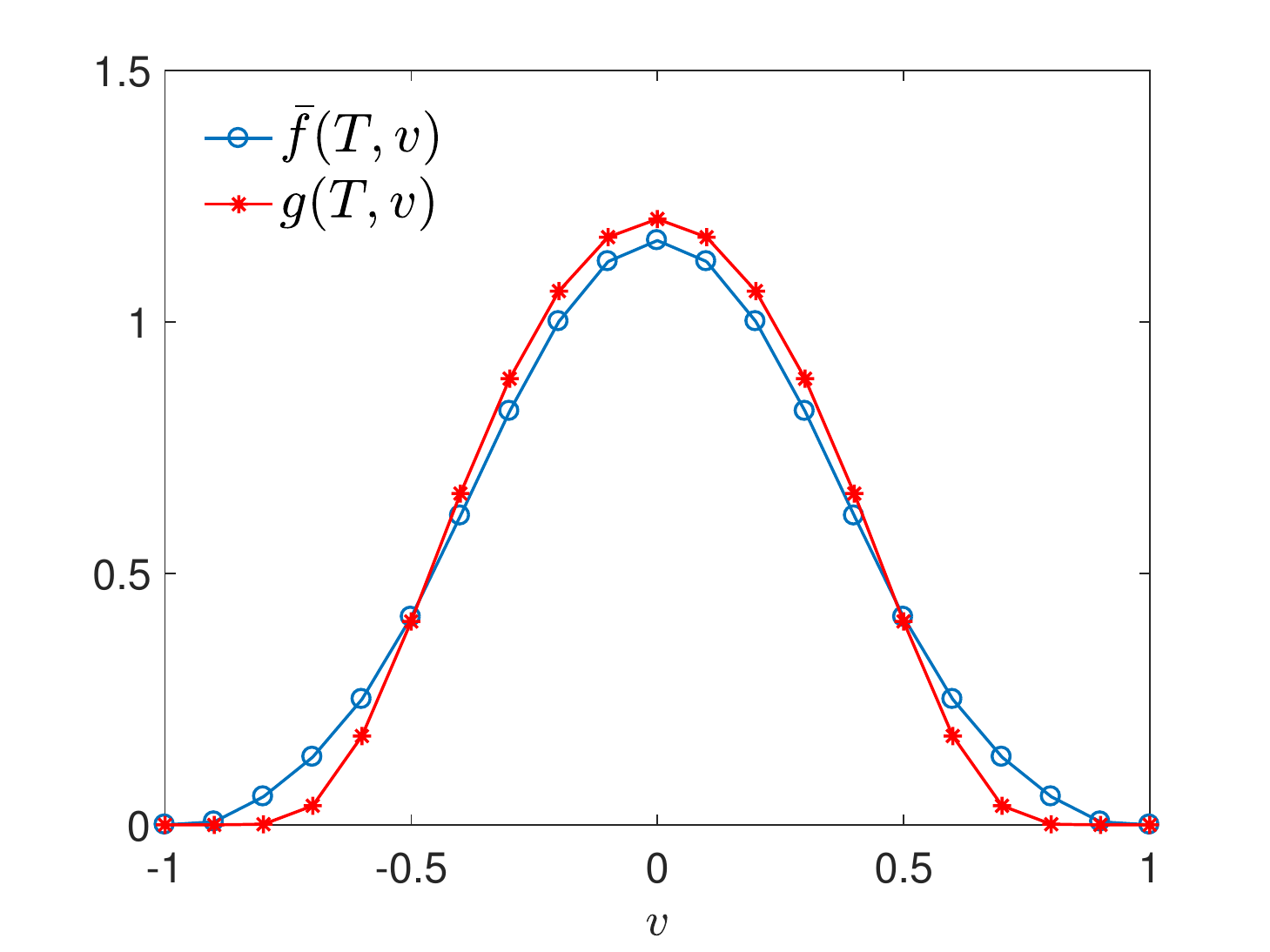}}
\subfigure[]{\includegraphics[scale=0.5]{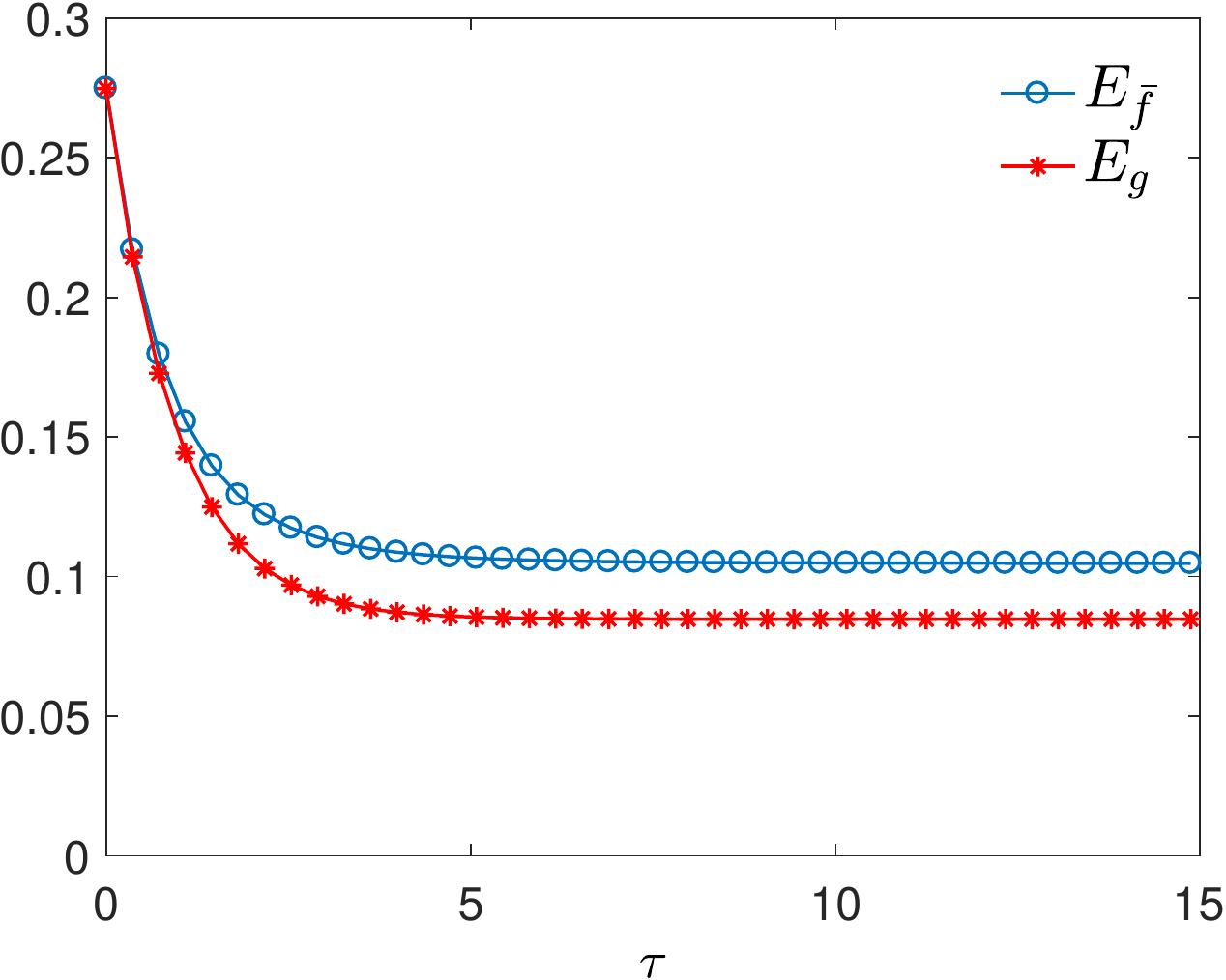}}
\caption{\textit{Example 2b}. (a) Asymptotic distributions of the deterministic and stochastic Fokker-Planck equations in the case of a consensus dynamics. (b) Evolution for $\tau\in [0,\,15]$ of the energies of the two models. In both cases the initial condition is the function $h_0$ in~\eqref{eq:h0} while the domain $[-1,\,1]$ is discretised with $N=21$ grid points and the time step $\Delta{\tau}=\Delta{v}^2/4=6.25\cdot 10^{-4}$ is used. We considered $q(\theta)$ like in~\eqref{eq:q} such that $q_0=0.5$ and $\lambda=0.4$. The evolutions of $\{f(\tau,\,v;\,\theta_k)\}_{k=0}^{10}$ and of $g(\tau,\,v)$ have been computed with $\text{SP}_\text{G}$ stochastic collocation methods.}
\label{fig:opinion_SS}
\end{figure}

\begin{figure}[!t]
\centering
\subfigure[]{\includegraphics[scale=0.5]{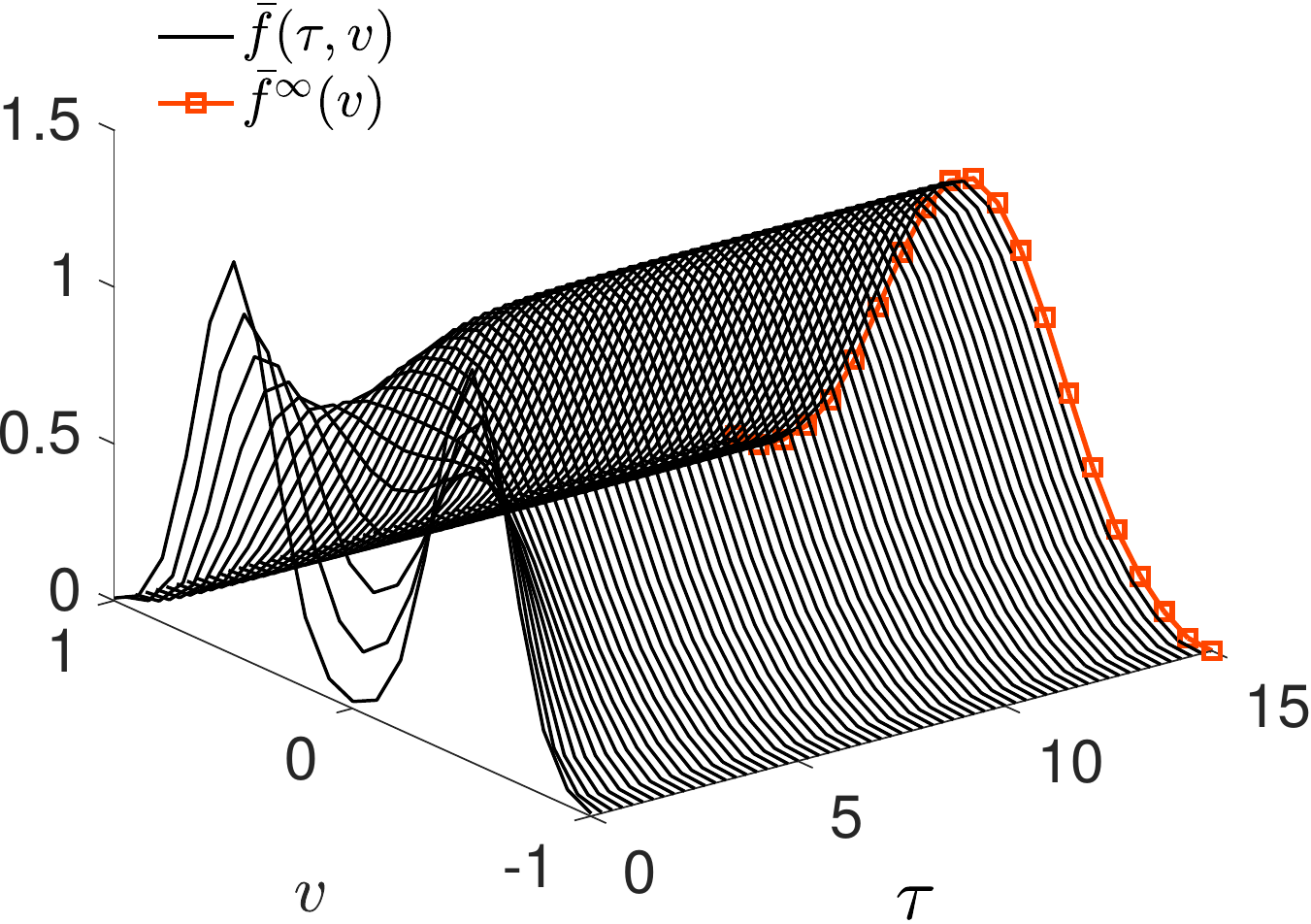}}
\subfigure[]{\includegraphics[scale=0.5]{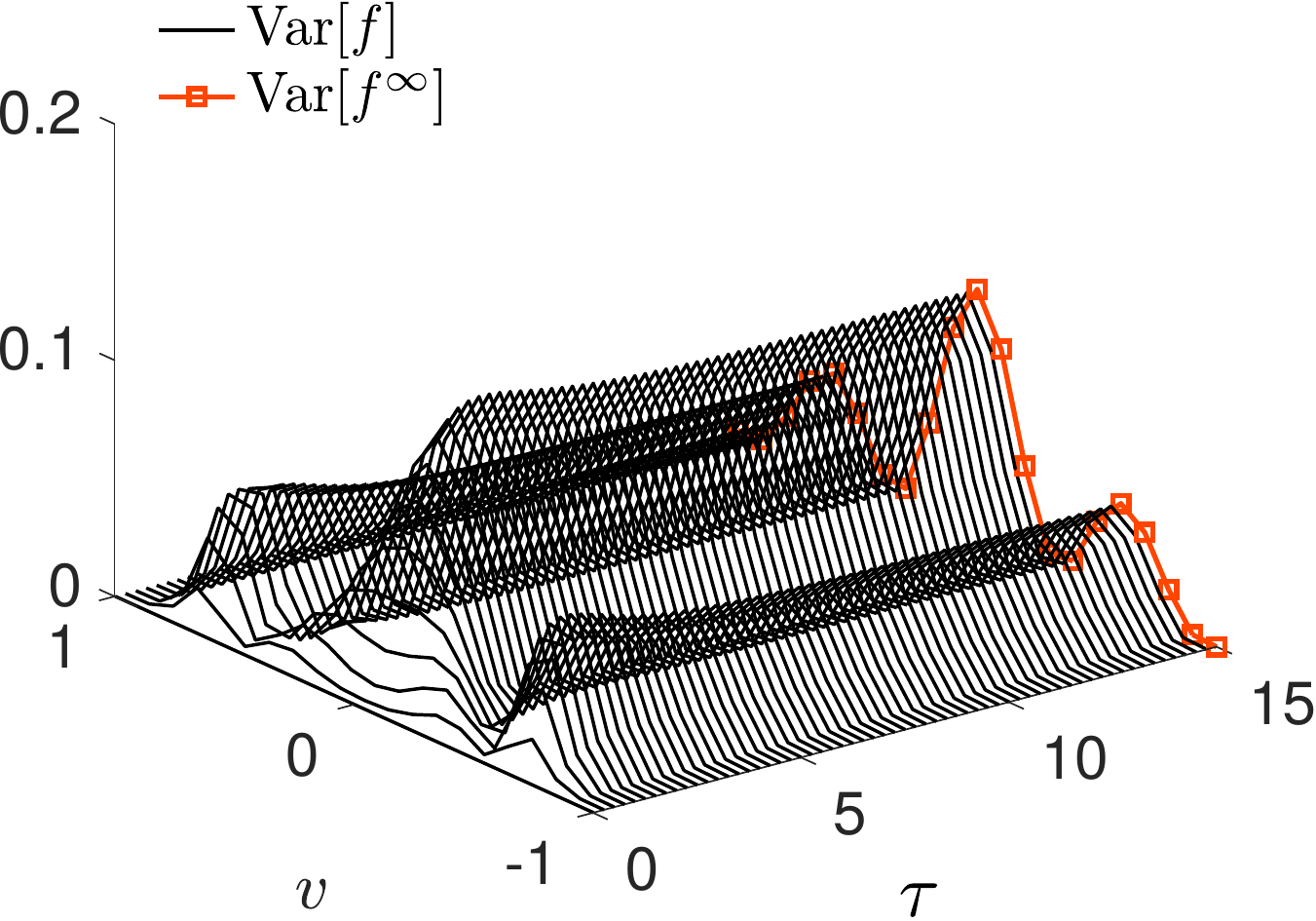}}
\caption{\textit{Example 2b}. Evolution of (a) $\bar{f}(\tau,\,v)$ and (b) $\operatorname{Var}[f](\tau,\,v)$ for the consensus model obtained with $M=10$ collocation points and the $\text{SP}_\text{G}$ scheme in the time interval $[0,\,15]$.}
\label{fig:opinion_surf}
\end{figure}

\subsubsection{Example 2b: Consensus with non-linear diffusion}
Now we consider the case discussed in Section~\ref{sect:non-linear_diff}, for which fully explicit analytical expressions of the functions $g^\infty$, $f^\infty(\cdot;\,\theta)$, $\bar{f}^\infty$ are not available in some cases. Structure preserving methods are therefore necessary in order to capture the correct large time solutions. Again we take $\theta\sim\mathcal{U}(-1,\,1)$ and we fix in particular $D(v)=D_0(1-v^2)$ with $D_0^2=0.025$. As initial condition we prescribe the deterministic distribution~\eqref{eq:h0} in the domain $V=[-1,\,1]$.

In Figure~\ref{fig:opinion_SS}a we show the functions $g(T,\,v)$ and $\bar{f}(T,\,v)$ for $T=15$, which have been obtained from the numerical solution of the deterministic and stochastic Fokker-Planck equations~\eqref{eq:FP_g},~\eqref{eq:FP_f} by means of the $\text{SP}_\text{G}$ scheme. The chosen value of the final time $T$ is such that they can be considered good approximations of the asymptotic distributions $g^\infty$ and $\bar{f}^\infty$, cf.~\eqref{eq:gf_non-lin_diff}. Furthermore, in Figure~\ref{fig:opinion_SS}b we show the evolution of the energies $E_{\bar f}$ and $E_g$ in the time interval $[0,\,T]$, $T=15$. Similarly to Figure~\ref{fig:FP_SS}b, we observe that the two models reach two different asymptotic energy values.

Finally, in Figure~\ref{fig:opinion_surf} we show the time evolution of the expected distribution function $\bar{f}(\tau,\,v)$ and of its $\theta$-variance, namely
\begin{align*}
	\operatorname{Var}[f](\tau,\,v) &:= \int_{\rangeth}f^2(\tau,\,v;\,\theta)\rho(\theta)\,d\theta-\left(\bar{f}(\tau,\,v)\right)^2 \\
	&= \frac{1}{2}\int_{-1}^{1}f^2(\tau,\,v;\,\theta)\,d\theta-\left(\bar{f}(\tau,\,v)\right)^2,
\end{align*}
of the stochastic Fokker-Planck equation~\eqref{eq:FP_consensus}$_2$ for $\tau\in [0,\,T]$. In particular, we report in red their asymptotic profiles at $\tau=T$.

\subsubsection{Example 3b: Bounded confidence model}
Among the most popular models of opinion dynamics in multi-agent systems, the \emph{bounded confidence model} has been first introduced in~\cite{hegselmann2002JASSS} and further investigated in~\cite{albi2016CHAPTER,albi2017KRM}. It describes a \emph{compromise} process in which the individuals interact only if their pre-interaction opinions $v,\,w\in V=[-1,\,1]$ are closer than a certain threshold $\Delta>0$, i.e. only if $\abs{w-v}\leq\Delta$. We are therefore again in the setting of Section~\ref{sect:non-linear_diff}, however with \emph{non-linear} binary interactions of the form~\eqref{eq:binary_symmetric_diff} with
$$ q=q(v,\,w;\,\theta)=\mathbbm{1}_{[0,\,\Delta(\theta)]}(\abs{w-v}). $$
Here $\Delta(\theta):=\Delta_0+\frac{1}{2}\theta$ is the bounded confidence threshold depending on the random input $\theta\sim\mathcal{U}(-1,\,1)$, with moreover $\frac{1}{2}\leq\Delta_0\leq\frac{3}{2}$ so as to guarantee $0\leq\Delta(\theta)\leq 2$ for all $\theta\in [-1,\,1]$, $\Delta=2$ being the maximum possible distance between any two opinions in $V$. In the following we fix specifically $\Delta_0=1$. Furthermore we notice that
$$ P(v,\,w;\,\theta)=q(v,\,w;\,\theta)(w-v) $$
in~\eqref{eq:P} and in particular
\begin{align*}
	\ave{P(v,\,w;\,\theta)} &= \frac{1}{2}\int_{-1}^{1}\mathbbm{1}_{[0,\,\Delta(\theta)]}(\abs{w-v})(w-v)\,d\theta \\
	&= \max\left\{\min\left\{1,\,\frac{3}{2}-\abs{w-v}\right\},\,0\right\}(w-v)
\end{align*}
in~\eqref{eq:FP_g} for the chosen value of $\Delta_0$.

Since analytical steady solutions are not available, we take advantage of SP numerical methods for a correct description of the large time trend of the distribution functions $g$ and $\bar{f}$.

\begin{figure}[!t]
\centering
\subfigure[]{\includegraphics[scale=0.5]{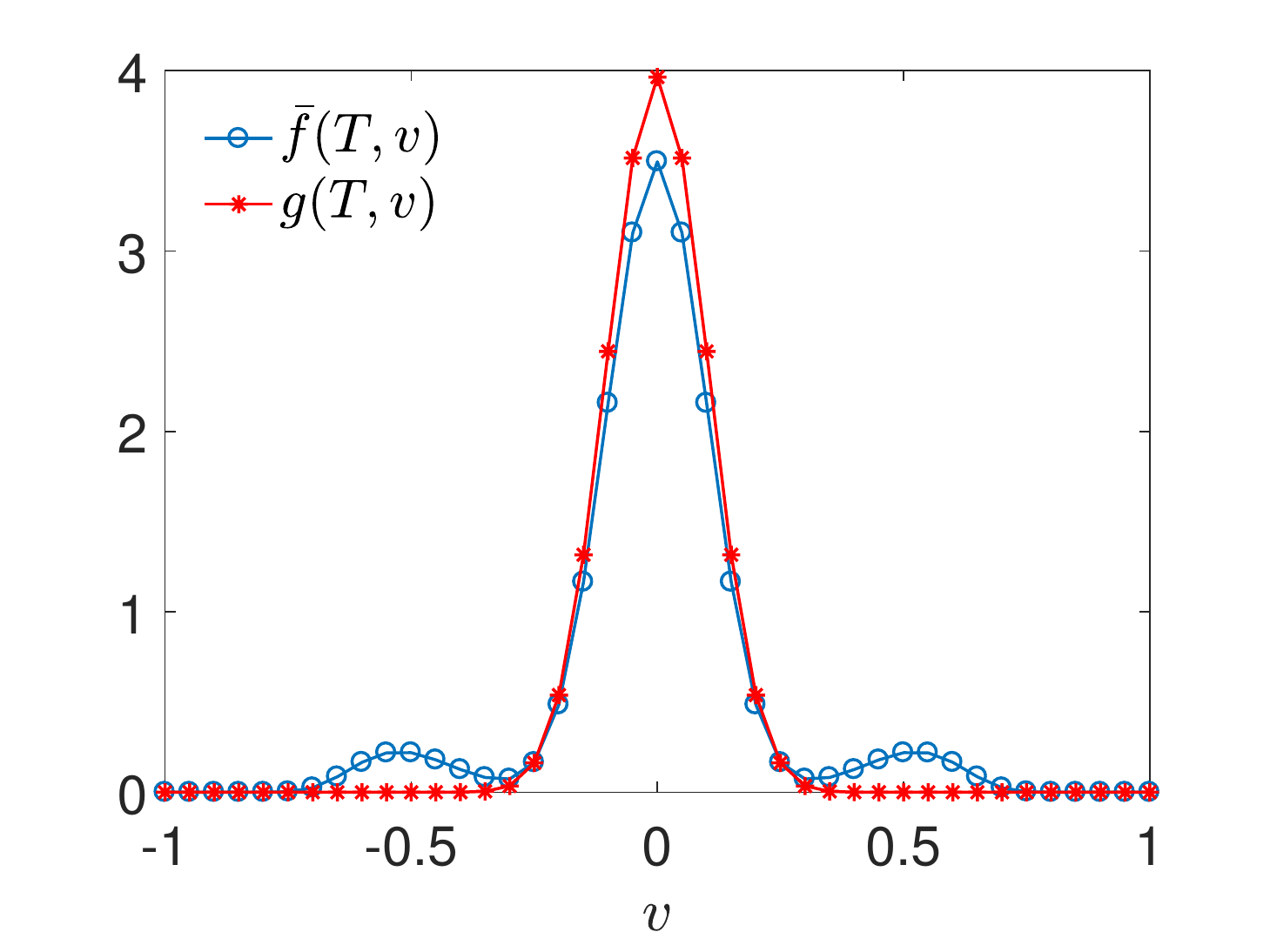}}
\subfigure[]{\includegraphics[scale=0.5]{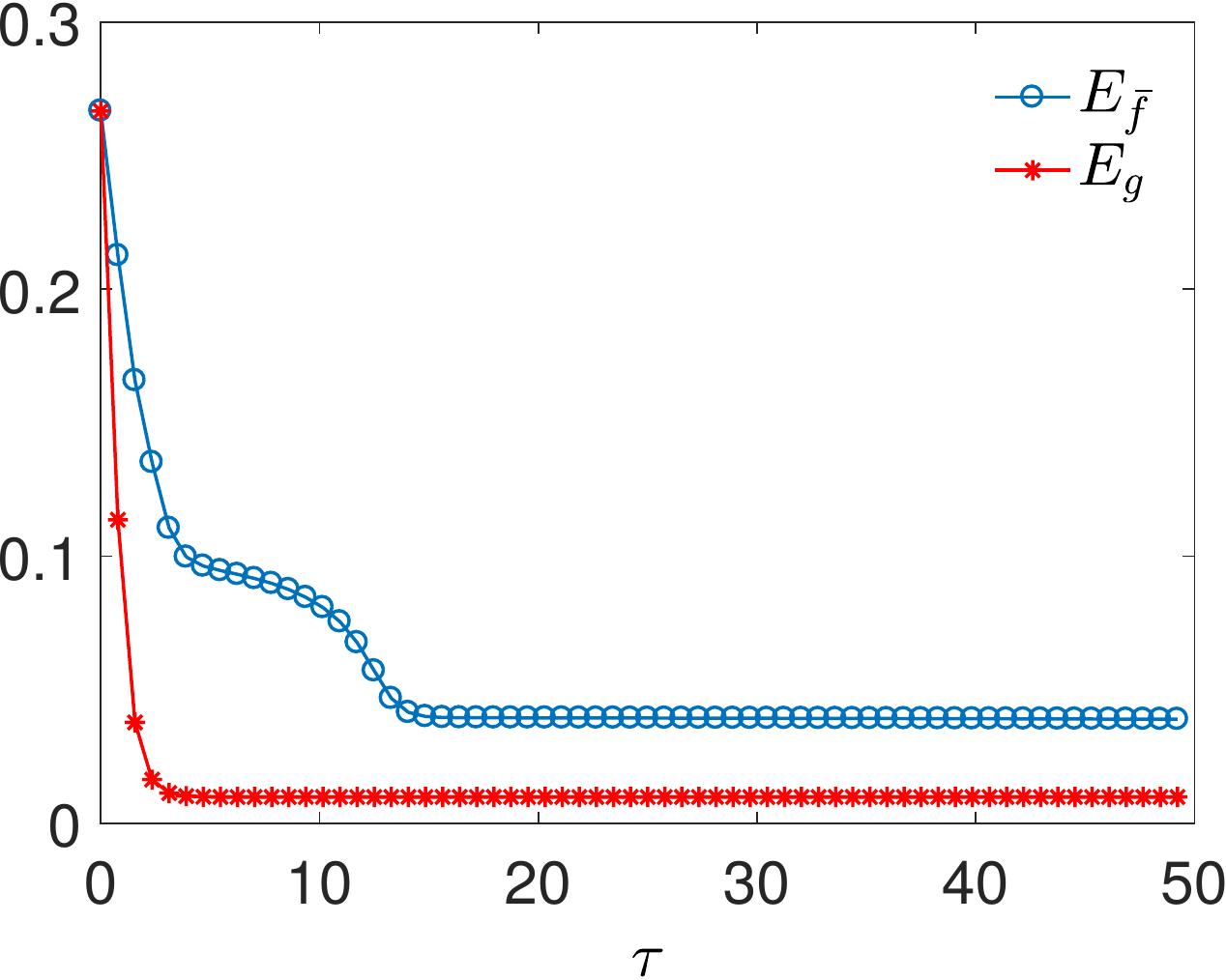}}
\caption{\textit{Example 3b}. (a) Asymptotic distributions of the deterministic and stochastic Fokker-Planck equations for the bounded confidence model. (b) Evolution in the time interval $[0,\,50]$ of the energies of the two models. In both cases the domain $V=[-1,\,1]$ is discretised with $N=21$ grid points and the time step $\Delta{\tau}=\Delta{v}^2/4\sigma^2$ with $\sigma^2=0.01$ is used. The evolutions of $\{f(\tau,\,v;\,\theta_k)\}_{k=0}^{10}$ and of $g(\tau,\,v)$ have been computed with $\text{SP}_\text{G}$ stochastic collocation methods.}
\label{fig:BC_A}
\end{figure}

In Figure~\ref{fig:BC_A}a we present a numerical approximation of $g^\infty$ and $\bar{f}^\infty$ computed by means of an $\text{SP}_\text{G}$ method in the time interval $[0,\,T]$, $T=50$, with $\Delta{\tau}=O(\Delta{v}^2)$ (see the caption of the figure for more details). The time integration has been performed by means of the fourth order Runge-Kutta method starting from the initial distribution $h_0$ in~\eqref{eq:h0}. Interestingly, the average asymptotic distribution of the stochastic model differs quite consistently from that of the deterministic model. In particular, besides the clustering about $v=0$, it shows two further clusters near $v=\pm\frac{1}{2}$, not displayed by the deterministic model, indicating that for some values of $\theta$ radicalisations of the opinions are possible.

In Figure~\ref{fig:BC_A}b we show the time evolution of the energies $E_g$, $E_{\bar{f}}$, which indeed tend asymptotically to two rather different values with quite different initial trends.

\begin{figure}[!t]
\centering
\subfigure[]{\includegraphics[scale=0.40]{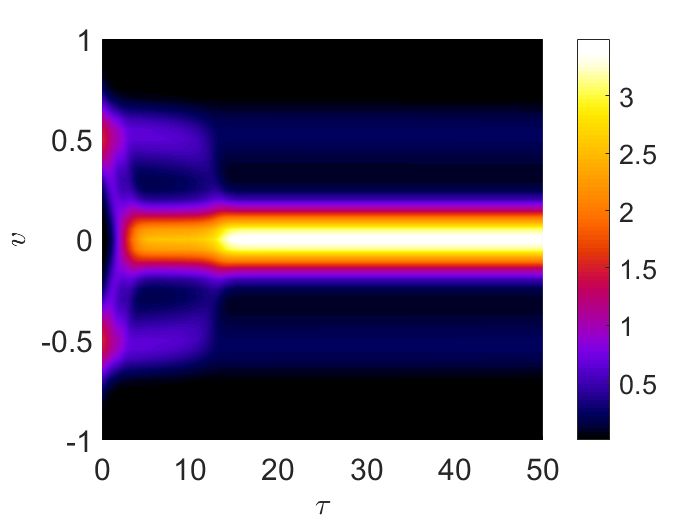}}
\subfigure[]{\includegraphics[scale=0.40]{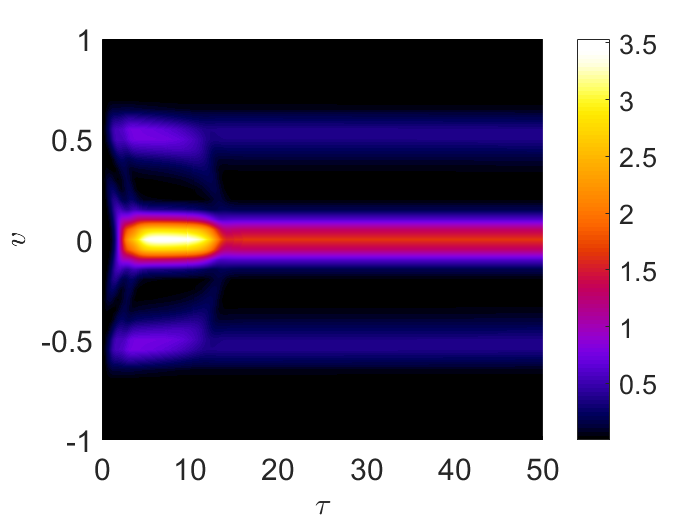}}
\caption{\textit{Example 3b}. Bounded confidence model: evolution for $\tau\in [0,\,50]$ of (a) $\bar{f}(\tau,\,v)$ and (b) $\operatorname{Var}[f](\tau,\,v)$ towards their asymptotic profiles obtained by means of the $\text{SP}_\text{G}$ scheme with $N=21$ grid points in the domain $V=[-1,\,1]$ and $M=10$ collocation points in $\rangeth=[-1,\,1]$. In this simulation the bounded confidence threshold is set to $\Delta(\theta)=\frac{3}{4}+\frac{1}{4}\theta$.}
\label{fig:BC}
\end{figure}

Finally, in Figure~\ref{fig:BC} we present the evolution of the expected distribution function $\bar{f}(\tau,\,v)$ and of its $\theta$-variance for the stochastic model in the time interval $[0,\,50]$ for a different choice of the parameters contained in $\Delta(\theta)$ (see the caption of the figure for details).

\section{Conclusions}
\label{sect:conclusions}
In this paper we have investigated the impact of uncertain microscopic interactions on the collective dynamics of multi-agent systems described by Boltzmann-type kinetic equations and their Fokker-Planck asymptotic limits. In the typical situation, the uncertainty is due to some parameters characterising the interactive dynamics, which are not known deterministically and which, in some cases, may drive the system towards quite different states.

We have proposed two different ways of quantifying the effect of such an uncertainty on the large-scale predictions of the models, depending on how the uncertainty is supposed to be distributed in the system. A first case is when the uncertain parameters vary from pair to pair of interacting particles, so that a common (although deterministically unknown) value of them valid for all the particles does not exist. Then one may average the binary interactions with respect to the probability distributions of the uncertain parameters and consider, at the kinetic level, the evolution of the distribution function subject to the representative \emph{average} binary interaction. A second case is instead when the uncertain parameters are in principle constant for all the particles. Then, at the kinetic level, one may consider the evolution of the distribution function \emph{parametrised} by the uncertain parameters. In this way one obtains a family of possible dynamics, which are more or less probable depending on the probability distributions of the parameters and which can be post-processed in order to extract statistical information such as the average trend and its variability. It is worth noticing that this second point of view is actually the one taken in uncertainty quantification (UQ).

By means of analytical and numerical methods we have shown that these two approaches are in general not equivalent and can lead to different interpretations of the collective dynamics, although in some cases, under suitable assumptions (e.g., in the absence of diffusive stochastic fluctuations), they may define the same large time trend of the system. Just to mention a relevant example, in case of general aggregation/consensus models they imply different conditions (generally more restrictive in the UQ-like setting) under which aggregation/consensus can be expected at the collective level in the two modelling approaches.

This demonstrates that, besides the quantification of the impact of the uncertainty, the interpretation of its role is first of all a source of complexity in multi-agent systems and a key issue in their mathematical modelling.

\section*{Acknowledgements}
A.T. is member of GNFM (Gruppo Nazionale per la Fisica Matematica) of INdAM (Istituto Nazionale di Alta Matematica), Italy. M.Z. is member of GNCS (Gruppo Nazionale per il Calcolo Scientifico) of INdAM, Italy.

The research that led to the present paper was partially supported by the research grant \emph{Numerical methods for uncertainty quantification in hyperbolic and kinetic equations} of GNCS-INdAM.

M.Z. acknowledges support from ``Compagnia di San Paolo'' (Torino, Italy).

\bibliographystyle{plain}
\bibliography{TaZm-boltzmannUQ}

\end{document}